\documentclass[11pt,twoside,a4paper,reqno]{amsart}

\usepackage[T1]{fontenc}
\usepackage{amsmath,amssymb,amsthm,mathtools}
\usepackage[expansion=false]{microtype}
\usepackage{aliascnt}
\usepackage[shortlabels]{enumitem}
\usepackage{booktabs,array}
\usepackage{xcolor}
\usepackage[colorlinks=true,linkcolor=blue!55!black,citecolor=blue!55!black,
            urlcolor=blue!55!black]{hyperref}
\usepackage[capitalize,noabbrev]{cleveref}

\newtheorem{theorem}{Theorem}[section]
\newaliascnt{proposition}{theorem}
\newtheorem{proposition}[proposition]{Proposition}
\aliascntresetthe{proposition}
\newaliascnt{lemma}{theorem}
\newtheorem{lemma}[lemma]{Lemma}
\aliascntresetthe{lemma}
\newaliascnt{corollary}{theorem}
\newtheorem{corollary}[corollary]{Corollary}
\aliascntresetthe{corollary}
\newaliascnt{conjecture}{theorem}

\aliascntresetthe{conjecture}
\theoremstyle{definition}
\newaliascnt{definition}{theorem}

\aliascntresetthe{definition}
\newaliascnt{problem}{theorem}

\aliascntresetthe{problem}
\theoremstyle{remark}
\newaliascnt{remark}{theorem}
\newtheorem{remark}[remark]{Remark}
\aliascntresetthe{remark}
\newaliascnt{example}{theorem}
\newtheorem{example}[example]{Example}
\aliascntresetthe{example}

\newcommand{\C}{\mathbb C}
\newcommand{\Z}{\mathbb Z}
\newcommand{\Q}{\mathbb Q}
\newcommand{\Pone}{\mathbb P^1}
\newcommand{\A}{\mathbb A}
\newcommand{\FIO}{\mathcal{FIO}}
\newcommand{\Mon}{\operatorname{Mon}}
\newcommand{\wt}{\operatorname{wt}}
\newcommand{\cD}{\mathcal D}
\newcommand{\cC}{\mathcal C}
\newcommand{\cE}{\mathcal E}

\newcommand{\Sing}{\operatorname{Sing}}
\newcommand{\PGL}{\operatorname{PGL}}

\newcommand{\Gm}{\mathbb G_m}
\newcommand{\Ga}{\mathbb G_a}
\newcommand{\rank}{\operatorname{rank}}
\newcommand{\Aut}{\operatorname{Aut}}
\newcommand{\Gr}{\operatorname{Gr}}

\newcommand{\Stir}[2]{S(#1,#2)}
\newcommand{\Sym}{\operatorname{Sym}}

\newcommand{\cS}{\mathcal S}
\newcommand{\cW}{\mathcal W}
\newcommand{\cT}{\mathcal T}

\numberwithin{equation}{section}
\hypersetup{
  pdftitle={Composition loci, invariant Bezoutians, and finite-orbit correspondences},
  pdfauthor={Ragni Piene and Boris Shapiro},
  pdfsubject={Algebraic correspondences and composition loci},
  pdfkeywords={finite orbit, Bezoutian, monodromy, polynomial decomposition}
}

\begin{document}

\title[Composition loci and invariant Bezoutians]
{Composition loci, invariant Bezoutians,\\ and finite-orbit correspondences}

\author[R.~Piene]{Ragni Piene}
\address{Department of Mathematics, University of Oslo, P.O. Box 1053 Blindern,
NO-0316 Oslo, Norway}
\email{ragnip@math.uio.no}

\author[B.~Shapiro]{Boris Shapiro}
\address{Department of Mathematics, Stockholm University, SE-106 91 Stockholm,
Sweden}
\address{Guangdong Technion--Israel Institute of Technology,
Shantou, Guangdong 515063, China}
\email{shapiro@math.su.se}

\date{\today}


\begin{abstract}
We compute projective degrees of composition loci for polynomial and rational
maps of the projective line. For normalized polynomials of degree $rs$ with
a right composition factor of degree $s$, a projective-bundle construction
gives the degree $S(r+s-2,r)$. We extend this calculation to composition
chains, and obtain $S(r+s-1,r)$ in the uncentered normalization used for
decomposable interpolation. For rational maps, a Grassmann-bundle model
gives a Pl\"ucker degree formula and an interpretation by prescribed critical
points. We relate these loci to factorization of the Bezoutian in the ring
of symmetric polynomials, where factors correspond to monodromy orbits on
unordered pairs. Fried's polynomial criterion then determines the invariant
reducibility locus. For rational maps we determine its composition components
when a marked pole fiber has exactly one ramified point and at least one
unramified point. The qualitative equivalence between finite full orbits and a
first integral is already known.
We give a direct proof from forward stabilization, with an explicit test
using bounded forward orbits and with quotient fibers identified at every
point. We also determine when absolute and invariant irreducibility differ
and give a growth formula for a class of translation correspondences.

\end{abstract}

\maketitle

\section{Introduction}\label{sec:intro}

An algebraic correspondence on the projective line is a rule that allows
a move from $z$ to $w$ whenever the pair $(z,w)$ lies on a specified algebraic
curve. A point may have several possible successors. Its full orbit consists
of all points reachable by repeating these moves in either direction.
For example, the rule $z^3=w^3$ confines every orbit to at most three points.
More generally, the relation $f(z)=f(w)$ for a nonconstant rational function
$f$ has the fibers of $f$ as its orbits. This raises a converse question.
If all full orbits of a correspondence are finite, are they the fibers of
some rational function?

We bring together results on finite-orbit correspondences and their
description by rational functions. We show that if enough starting points
have forward orbits of at most a prescribed size, then every full orbit is
finite. We give an explicit threshold for the number of starting points
and construct a rational function whose fibers agree with the full orbits
at every point.

Clearing denominators in $f(z)=f(w)$ and dividing out the diagonal factor
$z-w$ gives a polynomial called the Bezoutian.

The permutations of the points in a moving fiber determine the irreducible
factors of the Bezoutian and tell us which components have the whole fiber
as their generic orbit. In a decomposition $f=g\circ h$ with both degrees
greater than one, the Bezoutian of $h$ divides that of $f$, reflecting the
implication from $h(z)=h(w)$ to $f(z)=f(w)$.

This link motivates our main enumerative results. We compute the
projective degrees of the loci of polynomial and rational maps admitting
decompositions of prescribed degrees. The formulas give an exact count in
decomposable polynomial interpolation. They also count, with multiplicity,
rational maps with prescribed general critical points, up to changes of
the target coordinate. We describe the reducibility locus for rational
maps when a marked fiber has exactly one ramified point and at least one
unramified point.

The appendix asks how the number of reachable points grows when an orbit
is infinite. For a class of correspondences built from translations, we
show that the number of points reachable within $j$ steps grows on the
order of a power of $j$. The exponent is the rank of the group generated
by the translation steps.

We work over $\C$ unless another ground field is specified.

\subsection{Composition loci}
For a rational map $f=P/Q$ written in coprime form, the Bezoutian is
\[
 B_f(z,w)=\frac{P(z)Q(w)-P(w)Q(z)}{z-w}.
\]
Since it is symmetric in $z$ and $w$, it can be written uniquely as
$\overline B_f(u,v)$ with $u=z+w$ and $v=zw$. We write $B_P=B_{P,1}$ for
the polynomial case.

Two polynomials $P$ and $\widetilde P$ are \emph{affinely equivalent} if
$\widetilde P=L_1\circ P\circ L_2$ for polynomials $L_1,L_2$ of degree one.
Two rational maps are \emph{M\"obius equivalent} if they differ by M\"obius
transformations of the source and of the target. In a decomposition
$f=g\circ h$ with $\deg g,\deg h>1$, we call $g$ the \emph{left composition
factor} and $h$ the \emph{right composition factor}.

For a polynomial $P$ of degree $k\ge3$, Fried's theorem
\cite[Thm.~1, p.~45]{FriedSchur} and the Dickson factorization give a
complete criterion for reducibility in the symmetric coordinates
(\cref{th:fried}). The polynomial $\overline B_P$ is reducible precisely
when $P$ is decomposable, or when $k\ge5$ is prime and $P$ is affinely
equivalent to a Dickson polynomial $D_k(t,a)$. To study the geometry of
these loci, we work in the \emph{normalized slice}
\[
 \widetilde{\mathcal P}_k
 =\{t^k+a_{k-2}t^{k-2}+\cdots+a_1t\}\simeq\A^{k-2}
\]
of monic, centered polynomials with zero constant term. Every degree-$k$
polynomial is affinely equivalent to one in this slice. For $k=rs$ with
$r,s>1$, let $\cD_{r,s}\subset\widetilde{\mathcal P}_k$ denote the locus
of polynomials admitting a decomposition $G\circ H$ with $\deg G=r$ and
$\deg H=s$.

For a vector bundle $\cE$, we let $\mathbb P(\cE)$ denote the bundle of
one-dimensional subspaces, with tautological line bundle
$\mathcal O_{\mathbb P(\cE)}(-1)$, and $\Gr_2(\cE)$ the bundle of
two-dimensional subspaces; $\Gr(2,V)$ is the Grassmannian of two-dimensional
subspaces of a vector space $V$. We write $V_d=H^0(\Pone,\mathcal O(d))$ for
the space of binary forms of degree $d$.

\begin{theorem}[Composition geometry]\label{th:intro-scroll}
Let $r,s\ge2$.
\begin{enumerate}[{\rm (a)}]
\item For polynomials, the known normalized composition map
\cite[Theorem~2.2]{vzGMat2017} is a closed embedding
$\A^{r-1}\times\A^{s-2}\to\widetilde{\mathcal P}_{rs}$ with image $\cD_{r,s}$.
The closure of $\cD_{r,s}$ in $\mathbb P^{rs-2}\supset\widetilde{\mathcal P}_{rs}$
is the birational image of $\mathbb P(\cE)$, where
\[
 \cE=\mathcal O(-1)\oplus\cdots\oplus\mathcal O(-r)
 \quad\text{on }\mathbb P^{s-2},
\]
and it has degree $S(r+s-2,r)$, a Stirling number of the second kind.
\item For rational maps, the closure in $\Gr(2,V_{rs})$ of the locus of
pencils having a right composition factor of degree $s$ is the birational
image of the Grassmann bundle $\Gr_2(\Sym^r\cS)$ over $\Gr(2,V_s)$, where
$\cS$ is the tautological subbundle on $\Gr(2,V_s)$. It has codimension
$2(r-1)(s-1)$, and its Pl\"ucker degree is given by \eqref{eq:rational-degree}.
\end{enumerate}
\end{theorem}

For an ordered factorization $\mathbf d=(d_1,\ldots,d_\nu)$, we write
$\cC_{\mathbf d}$ for the locus of normalized compositions
$F_1\circ\cdots\circ F_\nu$ with $\deg F_i=d_i$. Its projective degree is
the product of Stirling numbers in \cref{cor:chain-degree}. The uncentered
interpolation degree and a closed formula for rational maps with a quadratic
right factor are given in \cref{cor:interpolation-degree,cor:rational-quadratic}.

\subsection{Correspondences and their orbits}
We now make the orbit question precise. A correspondence is a reduced
curve $\Gamma\subset\Pone\times\Pone$. Throughout, we assume that both
projections $\pi_1,\pi_2\colon\Gamma\to\Pone$ are finite. Equivalently,
$\Gamma$ has no \emph{fiber components}, which are components of the form
$\{a\}\times\Pone$ or $\Pone\times\{a\}$. We call
\[
 (n,m)=\bigl(\deg(\pi_1\colon\Gamma\to\Pone_z),\ \deg(\pi_2\colon\Gamma\to\Pone_w)\bigr)
\]
the \emph{projection degrees} of $\Gamma$. If $\Gamma$ is given by an affine
equation $F(z,w)=0$, then $F$ has degree $n$ in $w$ and degree $m$ in $z$.

The \emph{forward orbit} of a point $x\in\Pone$ is the set of points that can be
reached from $x$ by finitely many forward steps $z\mapsto w$; the
\emph{backward orbit} is defined using backward steps $w\mapsto z$; and the
\emph{full orbit} consists of all points reachable by finite sequences of
forward and backward steps in any order. All three contain $x$ itself. The
full orbit of $x$ is the connected component of $x$ in the graph on
$\Pone$ whose edges are the pairs in $\Gamma$. We call $\Gamma$ a
\emph{finite-orbit correspondence}, or simply $\FIO$, if all its full orbits
are finite. A \emph{rational first integral} of $\Gamma$ is a nonconstant
rational function $h$ such that $h(z)=h(w)$ for all $(z,w)\in\Gamma$, that is,
$h\circ\pi_1=h\circ\pi_2$. A rational first integral is constant on every
full orbit, so its existence bounds the size of every orbit by $\deg h$.

Any union of components of $f(z)=f(w)$ has first integral $f$.
The equivalence between finite full orbits and the existence of a first
integral follows from Bella\"iche's theorem \cite[Theorem~2.2.1]{Bellaiche}.
The formulation below also gives the forward-orbit criterion and describes
the quotient fibers. We discuss the relation to earlier work in
\cref{sec:prior}.

\begin{theorem}[Finite-orbit quotient]\label{th:intro-fio}
For a reduced correspondence $\Gamma\subset\Pone\times\Pone$ without fiber
components, the following conditions are equivalent:
\begin{enumerate}[{\rm (i)}]
\item uncountably many points have finite forward orbit;
\item every forward orbit is finite;
\item every forward orbit and every backward orbit is finite;
\item the full orbits are finite and their sizes are uniformly bounded;
\item $\Gamma$ admits a rational first integral;
\item $\Gamma$ is a union of components of the fiber relation of a rational
function.
\end{enumerate}
If these conditions hold, then $n=m$, and the forward, backward and full
orbits of every point coincide. Moreover, there is a \emph{minimal quotient}
$f$, unique up to postcomposition with an element of $\PGL_2(\C)$, whose
fibers, as sets, are exactly the orbits of $\Gamma$. Its degree is the size of
a generic orbit, and every rational first integral of $\Gamma$ belongs to
$\C(f)$. Each irreducible component of $\Gamma$ satisfies the same
conclusions, with its own minimal quotient.
\end{theorem}

The next result is an effective version of the implication
(i)$\Rightarrow$(vi). It is valid over any algebraically closed field of
characteristic zero, and it also yields a short proof of
\cref{th:intro-fio}.

\begin{theorem}[Effective forward-orbit test]\label{th:intro-effective}
Let $\Gamma$ have projection degrees $(n,m)$, let $N\ge1$, and put
$M_N=\sum_{i=0}^{N}m^i$. If more than
\[
 C(n,m,N)=2nN M_N
\]
points have forward orbits of size at most $N$, then all orbits of $\Gamma$
are the fibers of a rational function of degree at most $N$. If $\Gamma$ is
symmetric, the bound $C(n,m,N)$ may be replaced by $2nN(nN+1)$.
\end{theorem}

The distinction between forward and full orbits is essential. For the graph
of $w=z^2$, the points with finite forward orbit are $0$, $\infty$ and the
roots of unity. They form an infinite, Zariski-dense set, but there is no
rational first integral, since the generic orbit is infinite. In
particular, ``uncountably many'' in (i) cannot be replaced by ``infinitely
many''. Note also that no bound on the orbit size can depend only on $(n,m)$:
M\"obius transformations of finite order already give arbitrarily large finite
orbits in projection degrees $(1,1)$.

\subsection{Bezoutians and monodromy}
When a regular value of $f$ travels around a loop avoiding the branch
values, the points in its fiber are permuted. These permutations form the
monodromy group $\Mon(f)$. Its orbits on ordered pairs of distinct points
are called \emph{nontrivial orbitals} and correspond to the irreducible
components of $B_f=0$. In the symmetric coordinates $u,v$, the irreducible
factors correspond to orbits on unordered pairs. Thus $B_f$ is irreducible
exactly when $\Mon(f)$ is transitive on ordered pairs of distinct points
($2$-transitive). For $\overline B_f$, irreducibility is equivalent to
transitivity on unordered pairs ($2$-homogeneity). See
\cref{prop:monodromy}.

In the opening example $f(t)=t^3$, the Bezoutian $z^2+zw+w^2$ splits into two
linear factors over $\C$.
In symmetric coordinates it is the irreducible polynomial $u^2-v$.

\subsection{Organization}
Section~\ref{sec:fio} constructs the quotient from forward orbits.
Sections~\ref{sec:monodromy}--\ref{sec:polynomial} connect the quotient and
its Bezoutian to monodromy and composition. Section~\ref{sec:scroll}
computes the degrees of the composition loci. The following sections
study their intersections, classify the cases where absolute and
invariant irreducibility differ, and give small-degree examples.
The appendix studies orbit growth in translation models.

\subsection{Earlier work and the contribution of this paper}\label{sec:prior}

\emph{Composition varieties: known structure, new degree calculations.}
Normalized polynomial decompositions, their uniqueness, and polynomial
recovery of the factors belong to the classical and algorithmic theory of
decomposition \cite{Ritt,vonzurGathen1990,vonzurGathen2015}.
Von zur Gathen and Matera describe the affine geometry of these composition
varieties in \cite[Theorem~2.2]{vzGMat2017}. Their interpolation paper
\cite[Corollary~2.5 and Lemma~3.5]{vzGMat2024} already exhibits the Stirling
number $S(d+e-1,d)$ as a lower bound for the generic interpolation degree.
Our projective-bundle calculation gives the exact degree in both the
centered and uncentered normalizations. It also gives the chain-degree
product and the Grassmann-bundle formula for rational compositions.

\emph{Factorization and monodromy: known results and their invariant form.}
Components of a fiber product correspond to monodromy orbits; the
irreducibility criterion for a residual fiber curve is stated in
\cite[Corollary~2.3]{Pak2011}. Passing from ordered to unordered pairs gives
the factorization in $\C[z+w,zw]$ used here. These coordinates also occur in
\cite[proof of Theorem~1.1]{MirmanShukla}. The polynomial classification is
Fried's theorem \cite[Theorem~1]{FriedSchur}, followed by the classical
Dickson difference identity \cite[Proposition~1.7]{Turnwald}. Ritt--Engstr\"om
theory describes intersections of decomposition loci
\cite{Engstrom,Tortrat,ZieveMueller,vonzurGathen,Ziegler}; we retain only
the consequences needed for their degrees and invariant factors.

\emph{Rational maps: consequences of standard group theory and a component calculation.}
The simple-pole criterion uses Jordan's theorem as recalled in
\cite[\S2]{Jones}. Generic wreath-product monodromy follows from the
standard isolated-inertia argument, familiar from
\cite{Odoni,Juul,JuulKuMaTu,KonigNeRo}. Our contribution in this setting is
the explicit description and dimension of the composition components with
a prescribed pole profile. Kantor's theorem \cite{Kantor}, together with
the genus-zero condition, governs the discrepancy between absolute and
invariant irreducibility. The cyclic and elliptic maps that occur are
classical \cite{Ritt3,Burda,GuMuSa}. We combine Kantor's theorem with the
affine genus-zero classification \cite{Neubauer1992,MagaardShpectorovWang}
to exclude the nonprime semilinear possibilities.

\emph{Finite orbits: known qualitative theorem, a direct proof and an effective test.}
Bella\"iche proves the equivalence between a first integral and infinitely
many finite complete sets, as well as the uniform boundedness of all full
orbits \cite[Theorem~2.2.1]{Bellaiche}. For irreducible correspondences,
Bell--Moosa--Topaz provide a broader invariant-hypersurface theorem
\cite[Corollary~6.2]{Bell}. Fried's earlier theorem
\cite[Theorem~1]{FriedPoncelet} constructs a common quotient for finite
alternating correspondences; his exceptional points are discussed in
\cite[\S1.B]{FriedPoncelet}. A common quotient may use different functions
on the two coordinates and need not be an iteration first integral.
For example, the graph $w=z^2$ has a common quotient but infinite generic
iteration orbits. Beaumont discusses the relation between forward and full
orbits in \cite[Remark~4.4]{Beaumont}; exchanging a correspondence with its
transpose, however, does not preserve a given forward-orbit hypothesis.
Our proof treats either orientation directly and identifies the quotient
fibers at every point. The bound is useful particularly when $n<m$.
For $n>m$, Bella\"iche's ramification bound is stronger; in the balanced
case his bounds concern complete classes, whereas ours counts points with
bounded forward orbit. The comparison is made explicitly in
\cref{rem:bound-comparison}.

\section{Forward stabilization and the quotient}\label{sec:fio}

In this section $\mathbb K$ denotes an algebraically closed field of
characteristic zero; assertions involving countability are made over $\C$.
Let $\Gamma \subset \Pone \times \Pone$ be a reduced curve without fiber
components, with projection degrees $(n,m)$. We compose relations
set-theoretically: for $R,S\subset\Pone\times\Pone$,
\[
 S\circ R=\{(x,z):\ (x,y)\in R\ \text{and}\ (y,z)\in S\ \text{for some } y\}.
\]
Set $\Gamma^0=\Delta$, the diagonal, and $\Gamma^{i+1}=\Gamma\circ\Gamma^i$ for
$i\ge0$, so that $\Gamma^i(x)$ is the set of points reachable from $x$ in
exactly $i$ forward steps. Put
\[
 F_j=\bigcup_{i=0}^{j}\Gamma^i,
 \qquad e_j=\deg(F_j\to\Pone_z),\qquad
 \hat e_j=\deg(F_j\to\Pone_w).
\]
Each $F_j$ is a reduced closed curve without fiber components. Indeed, the
space of paths of length $i$ in $\Gamma$ is an iterated fiber product of the
finite flat maps $\pi_1,\pi_2$ over the smooth curve $\Pone$; it is therefore
finite and flat over the initial point, hence of pure dimension one, and
$\Gamma^i$ is its image under the proper map recording the two endpoints. This
image is closed, and it is finite over either endpoint. In particular, it has no
fiber components and no isolated points. By definition, the fiber
$F_j(x)$ is the set of points reachable from $x$ in at most $j$ forward steps.

The elementary stabilization below is related to the path argument in
\cite[Lemma~4.4.6]{Bellaiche}. Fried's earlier construction
\cite[Theorem~1]{FriedPoncelet} uses alternating correspondences and a field
intersection; our relations retain the forward direction throughout.

\begin{lemma}[Stabilization]\label{lem:stab}
The chain $F_0\subseteq F_1\subseteq\cdots$ has the following properties:
\begin{enumerate}[{\rm (a)}]
\item a strict inclusion $F_j\subsetneq F_{j+1}$ increases each projection
degree by at least one;
\item $e_{j+1}\le1+ne_j$ and $\hat e_j\le\sum_{i=0}^jm^i$;
\item if $F_j=F_{j+1}$, then $F_i=F_j$ for all $i\ge j$, and $F_j$ equals the
forward relation $\bigcup_{i\ge0}\Gamma^i$.
\end{enumerate}
\end{lemma}
\begin{proof}

A component of $F_{j+1}$ not contained in $F_j$ has positive degree over both
coordinates, which gives (a). The identity
$F_{j+1}=\Delta\cup(\Gamma\circ F_j)$ yields $e_{j+1}\le 1+ne_j$, and
$\hat e_j\le\sum_{i\le j}m^i$ holds because $\Gamma^i$ has degree at most
$m^i$ over the second coordinate. The same identity shows that an equality
$F_j=F_{j+1}$ propagates to all later steps, proving (c).
\end{proof}

The next lemma is a useful algebraic form of the finite-set counting
principle in \cite[Scholium~2.1.2]{Bellaiche}. Its proof works away from a
finite exceptional set and then uses closedness to include that set.

\begin{lemma}[Finite algebraic preorders]\label{lem:quasi}
Let $X$ be a smooth projective curve and $E\subset X\times X$ a reduced
curve without fiber components. If $\Delta\subset E$ and $E\circ E=E$,
then $E$ is symmetric, hence an equivalence relation.
\end{lemma}
\begin{proof}

Let $e$ be the degree of the first projection $E\to X$. Outside a finite set
$Z\subset X$, every fiber $E(x)=\{y:(x,y)\in E\}$ has exactly $e$ elements.
The set $E^{-1}(Z)=\{x:E(x)\cap Z\ne\emptyset\}$ is also finite. Let
$x\notin Z\cup E^{-1}(Z)$ and $y\in E(x)$. Transitivity gives
$E(y)\subseteq E(x)$, and both sets have $e$ elements because $x,y\notin Z$;
hence $E(y)=E(x)$. Since $x\in E(x)$ by reflexivity, we get $(y,x)\in E$. The
pairs $(x,y)$ considered are dense in every component of $E$, because $E$ has
no fiber components. Since $E$ is closed, $E^{-1}\subseteq E$, and therefore
$E^{-1}=E$.
\end{proof}

For a relation $E$ we call a rational function $h$ \emph{$E$-invariant} if
$h(x)=h(y)$ for all $(x,y)\in E$. The $E$-invariant rational functions,
together with the constants, form a subfield of $\mathbb K(t)$.

We include the quotient construction to make clear why the fibers agree
with the equivalence classes also at exceptional points. Compare the
field-theoretic quotient in \cite[Theorem~1]{FriedPoncelet}; Bella\"iche's
Lemma~4.4.6 is a stabilization lemma, rather than a quotient theorem.

\begin{lemma}[Quotient]\label{lem:quotient}
Let $E\subset\Pone\times\Pone$ be a reduced equivalence relation without fiber
components, and let $e=\deg(E\to\Pone)$. Then $E$ is the reduced fiber
relation $\{(x,y):f(x)=f(y)\}$ of a rational function $f$ of degree $e$,
unique up to postcomposition with a M\"obius transformation, and the field of
$E$-invariant rational functions is $\mathbb K(f)$.
\end{lemma}
\begin{proof}

In an affine chart, the class $E(z)$ of a generic point is the zero set of a
monic polynomial $w^e+c_1(z)w^{e-1}+\cdots+c_e(z)$ with
$c_i\in\mathbb K(z)$. The coefficients $c_i$ are $E$-invariant, since
equivalent points have the same class, and they are not all constant, since a
generic class contains $z$ and therefore cannot be independent of $z$. Let
$L=\mathbb K(c_1,\ldots,c_e)$. Two general points have the same coefficient
vector exactly when they lie in the same class, so the map
$z\mapsto(c_1(z),\ldots,c_e(z))$ has degree $e$ onto its image and
$[\mathbb K(z):L]=e$. By L\"uroth's theorem, $L=\mathbb K(f)$ for a rational
function $f$ of degree $e$. Every $E$-invariant rational function is constant
on the fibers of the coefficient map, hence factors through it and lies in
$L$.

The equality $f(z)=f(w)$ holds at general points of every component of $E$,
hence everywhere on $E$, as an equality of morphisms to $\Pone$. Thus $E$ is
contained in the reduced fiber relation of $f$. Both curves have first
projection degree $e$ and neither has fiber components, so they coincide; in
particular, the classes of $E$ are exactly the fibers of $f$, also at the
finitely many exceptional points. Finally, two generators of the field
$L=\mathbb K(f)$ differ by a M\"obius transformation.
\end{proof}

\begin{proof}[Proof of \cref{th:intro-effective}]\label{subsec:effective}
Let $D$ be a reduced curve without fiber components and with projection
degrees $(e,\hat e)$. The discriminant, with respect to $w$, of a bihomogeneous equation of $D$ is a binary form in $z$ of degree $(2e-2)\hat e$, which is nonzero because the generic fiber of $D\to\Pone_z$ consists of $e$ distinct
points. Hence
\begin{equation}\label{eq:Zbound}
 \#\{x\in\Pone:\ |D(x)|<e\}\le(2e-2)\hat e.
\end{equation}
This includes the case $e=1$, in which every fiber consists of one point.

Suppose that $e_j>N$ for some $j$, and let $j$ be the smallest such index.
Then $j\le N$. Indeed, if $F_i=F_{i+1}$ for some $i<j$, then $e_j=e_i\le N$ by
\cref{lem:stab}(c), which is impossible; so all inclusions before step $j$ are
strict, and by \cref{lem:stab}(a) and $e_0=1$ we have $e_i\ge i+1$ for
$i\le j$, which forces $j\le N$. By \cref{lem:stab}(b) we obtain
$e_j\le 1+ne_{j-1}\le nN+1$ and $\hat e_j\le M_N$. Now let $x$ be a point whose
forward orbit has at most $N$ elements. Since $F_j(x)$ is contained in this
orbit, we have $|F_j(x)|\le N<e_j$, so $x$ lies in the exceptional set
\eqref{eq:Zbound} for $D=F_j$. That set has at most
$(2e_j-2)\hat e_j\le 2nNM_N$ elements, contradicting the hypothesis.

Thus $e_j\le N$ for all $j$. By \cref{lem:stab}(a), the chain becomes
stationary, and by \cref{lem:stab}(c) its final term is the forward relation
$E=\bigcup_{i\ge0}\Gamma^i$, whose fibers are the forward orbits. The relation
$E$ is reflexive and transitive, so it is an equivalence relation by
\cref{lem:quasi}, and by \cref{lem:quotient} its classes are the fibers of a
rational function of degree $e\le N$. Since $E$ is symmetric, forward orbits
are also backward orbits, and hence they are the full orbits.

If $\Gamma$ is symmetric, then so is every $F_j$, and $\hat e_j=e_j\le nN+1$.
At the first crossing, the bound \eqref{eq:Zbound} is then at most
$2nN(nN+1)$.
\end{proof}

\begin{proof}[Proof of \cref{th:intro-fio}]
The implications (ii)$\Rightarrow$(i) and (iii)$\Rightarrow$(ii) are trivial,
and (iv)$\Rightarrow$(iii) holds because forward and backward orbits are
contained in full orbits. A rational first integral $h$ is constant on full
orbits, so all full orbits have at most $\deg h$ elements; this proves
(v)$\Rightarrow$(iv). If $\Gamma$ is a union of components of the fiber
relation of $f$, then $f$ is a rational first integral, which proves
(vi)$\Rightarrow$(v).

Assume (i). The set of points with finite forward orbit is the union, over
$N\ge1$, of the sets of points whose forward orbit has at most $N$ elements.
If it is uncountable, one of these sets is infinite, and the proof of
\cref{th:intro-effective} shows that the forward relation $E$ is a closed
curve, that it is an equivalence relation, and that its classes are the fibers
of a rational function $f$. Since $\Gamma\subseteq E$, this proves (vi), and
we have also seen that forward, backward and full orbits coincide at every
point. The function $f$ is a rational first integral. Conversely, a rational first 
integral is constant on forward orbits, hence $E$-invariant, and by \cref{lem:quotient} 
it lies in $\C(f)$. A nonconstant element $g=h\circ f$ of $\C(f)$ has degree
$(\deg h)(\deg f)$, which proves the minimality and uniqueness of $f$.

Let $C$ be an irreducible component of $\Gamma$, and let $\pi_1,\pi_2$ be the
two projections of its normalization. Since $f\circ\pi_1=f\circ\pi_2$,
multiplicativity of degrees gives $\deg\pi_1=\deg\pi_2$; this shows $n=m$.
The orbits of $C$ are contained in those of $\Gamma$, hence bounded, so the
whole argument applies to $C$ as well.
\end{proof}

\begin{corollary}[Balance]\label{cor:balance}
Let $\Gamma$ be a finite-orbit correspondence. Every forward-invariant subset
of $\Pone$ is also backward invariant, and each irreducible component of
$\Gamma$ has equal projection degrees.
\end{corollary}

\begin{remark}[Comparison with earlier bounds]\label{rem:bound-comparison}
For $n>m$, the forward version of Bella\"iche's ramification estimate
\cite[Proposition~2.1.1 and its dual statement]{Bellaiche} is stronger.
If $R_1$ is the total ramification of the first projection from the
normalization of $\Gamma$, the number of points with finite forward orbit
is at most
\[
 \frac{R_1}{n-m}\le\frac{2m(n-1)}{n-m},
\]
independently of $N$. Indeed, a finite forward-invariant set $S$ satisfies
$n|S|-R_1\le |\pi_1^{-1}(S)|\le |\pi_2^{-1}(S)|\le m|S|$;
apply this to finite unions of finite forward orbits. The final estimate
uses the bidegree bound on the discriminant.
When $n=m$, Bella\"iche bounds the number of finite complete classes
\cite[\S3 and Theorem~4.5.3]{Bellaiche}, whereas our test bounds the number
of points whose forward orbits have prescribed bounded size. The two
bounds are therefore not directly comparable numerically. When $n<m$,
finite forward orbits may be infinite in number but countable (as for
$z\mapsto z^2$), so the one-sided bounded-size formulation is essential.
\end{remark}

\begin{remark}\label{rem:field}
The effective test is purely algebraic; for instance, it applies over
$\overline{\Q}$. The formulation with uncountably many points, in contrast,
requires an uncountable ground field, and we draw no conclusion from
infinitely many finite forward orbits of unbounded size over a countable
field. The hypothesis that $\Gamma$ has no fiber components cannot be dropped:
the relation $(w-c)(w-z)=0$, which contains the horizontal line $w=c$, has
finite forward orbits, but the backward orbit of $c$ is the whole line.
\end{remark}

\subsection{An arbitrary smooth projective curve}\label{subsec:anycurve}

Let $X$ be a smooth projective curve over $\C$ and $\Gamma \subset X\times X$
a reduced curve without fiber components. Forward, backward and full orbits are
defined as before.

\begin{proposition}\label{th:anycurve}
Conditions {\rm (i)--(iv)} of \cref{th:intro-fio} are equivalent for
$\Gamma\subset X\times X$, and they hold if and only if there is a
nonconstant morphism $f\colon X\to Y$ to a smooth projective curve $Y$ whose
fibers are the orbits of $\Gamma$. Such a morphism of minimal degree is unique
up to an isomorphism of its target, and the two projection degrees of every
component of $\Gamma$ are equal.
\end{proposition}
\begin{proof}
The construction of the path curves $F_j$ and \cref{lem:stab,lem:quasi} carry
over unchanged. If the degrees $e_j$ were unbounded, every point outside the
countable union of the finite exceptional sets of the curves $F_j$ would have
infinite forward orbit. Hence uncountably many finite forward orbits force the
chain to stabilize, and the stable relation $E$ is an equivalence relation.
Let $e$ be its first projection degree. Sending a general point $x$ to its
class $E(x)$ defines a nonconstant rational map $X\dashrightarrow\Sym^eX$;
let $Y$ be the normalization of its image. The resulting morphism
$f\colon X\to Y$ has degree $e$, and its generic fibers are the classes of
$E$. As in the proof of \cref{lem:quotient}, the reduced fiber relation of
$f$ and $E$ have the same projection degree and no fiber components, so they
are equal. The remaining assertions follow, as before, from the descent of
invariant functions and the multiplicativity of degrees.
\end{proof}

The proof of \cref{th:intro-effective} also gives an effective version on $X$.
The only additional ingredient is a bound for the number of fibers with
fewer than the expected number of points.

\begin{proposition}[An effective bound on $X$]\label{rem:curve-bound}
Let $X$ be a smooth projective curve of genus $g$ over an algebraically closed
field of characteristic zero, and let $\Gamma\subset X\times X$ be a reduced
curve without fiber components, with projection degrees $(n,m)$. For $N\ge1$
put $E_N=nN+1$, $M_N=\sum_{i=0}^Nm^i$ and
\[
 C_X:=\left\lceil\frac{(E_N+M_N)^2}{2}\right\rceil
       +\max(2g-2,0)\,M_N .
\]
If more than $C_X$ points of $X$ have forward orbits of size at most $N$, then
the forward relation of $\Gamma$ is an equivalence relation whose classes are
the fibers of a morphism $X\to Y$ of degree at most $N$.
\end{proposition}

\begin{proof}
The curve $D\subset X\times X$ is reduced with projection degrees $(e,\hat e)$ and has no fiber comoponents. Hence the map $\pi_{1}\colon D\to X$ is flat, and $\mathcal V:=\pi_{1*}\mathcal O_D$ is a vector bundle of rank $e$. We claim that the following holds:
\begin{multline*}
  \#\{x\in X:\ |D(x)| < e\} \le -2\deg \mathcal V=D^2+(2g-2)\hat e\\
 \le\frac{(e+\hat e)^2}{2}+\max(2g-2,0)\hat e.
\end{multline*}
The trace pairing $ \mathcal V\otimes \mathcal V\to\mathcal O_X$ gives a map $\mathcal V\to \mathcal V^\vee$. Its determinant $\delta$ is therefore a section of
$
 \det \mathcal V^\vee\otimes\det \mathcal V^\vee=(\det \mathcal V)^{-2}$, so 
 $ \deg\delta=-2\deg \mathcal V$.
Note that $\deg \mathcal V\le 0$, since $\deg\delta \ge 0$. Since a generic fiber of $\pi_1$ is reduced and we are in characteristic 0, $\pi_1$ is generically \'etale, so $\delta \ne 0$. The zeros of $\delta$ are the points $x\in X$ such that the fiber $D(x)$ is nonreduced, which happens precisely when $|D(x)|< e$. Indeed, multiplication by a nilpotent element has trace zero.
Hence
\[
 \#\{x\in X:\ |D(x)|<e\}\le-2\deg \mathcal V .
\]
By Riemann--Roch on $X$,
\[
 \chi(\mathcal O_D)=\chi(X,\mathcal V)=\deg \mathcal V+e(1-g).
\]
By adjunction on $X\times X$, with $K_{X\times X}=\pi_1^*K_X+\pi_2^*K_X$,
\[
 2p_a(D)-2=D^2+K_{X\times X}\cdot D=D^2+(2g-2)(e+\hat e).
\]
Since $\chi(\mathcal O_D)=1-p_a(D)$, combining the two equations gives
\[
 -2\deg \mathcal V=2p_a(D)-2+2e(g-1)=D^2+(2g-2)\hat e .
\]
Note that in the case $X=\Pone$, we have $D^2=2e\hat e$, hence we get 
\[
 \#\{x\in X:\ |D(x)|<2e\hat e -2\hat e=(2e-2)\hat e,
\]
which is (2.1).

For the last inequality, we apply the Hodge index theorem to the sum of the the two fiber classes, 
whose self-intersection is two. Apply this bound at the first index
$e_j>N$, where $e_j\le E_N$ and $\hat e_j\le M_N$, and repeat the proof
of \cref{th:intro-effective}. 
\end{proof}

We note that the branch divisor of the normalization of $D$ alone would not
suffice in this argument: it does not detect fibers in which two distinct
branches of $D$ pass through the same point.

\begin{remark}\label{rem:symmetric-complete}
For symmetric correspondences, infinitely many finite orbits already imply
the existence of a first integral by Bella\"iche's theorem
\cite[Theorem~2.2.1, p.~1878]{Bellaiche}, since for symmetric $\Gamma$ a finite forward
orbit is a finite complete set. This stronger qualitative statement is
independent of the explicit bounded-orbit test above.
\end{remark}

\section{Monodromy and the orbital dictionary}\label{sec:monodromy}

Let $f=P/Q\colon\Pone\to\Pone$ be a rational function of degree $N\ge2$, with
$P,Q$ coprime. In bihomogeneous coordinates, $P(z)Q(w)-P(w)Q(z)=0$ is the
fiber product of $f$ with itself. The diagonal is one of its components, and
the residual curve is the closure of $B_{P,Q}=0$; we call it the
\emph{residual fiber curve}. In characteristic zero this fiber product is generically reduced, hence
its equation and $B_{P,Q}$ are squarefree. The diagonal is not a component
of the residual curve, since $B_{P,Q}(w,w)=P'Q-PQ'$ is not identically zero.

The ordered-pair statement is standard; see
\cite[Proposition~2.1 and Corollary~2.3]{Pak2011}. We spell out its
unordered-pair consequence because our factorization ring is
$\C[z,w]^{\mathfrak S_2}=\C[z+w,zw]$. These invariant coordinates already
occur in \cite[proof of Theorem~1.1, p.~90]{MirmanShukla}.

\begin{proposition}[Residual fiber curves]\label{prop:monodromy}
The irreducible components of $B_{P,Q}=0$ correspond bijectively to the
nontrivial orbitals of $G:=\Mon(f)$ on a generic fiber. If an orbital has
size $Nd$, the corresponding component has projection degrees $(d,d)$. In
particular,
\[
 B_{P,Q}\ \text{is irreducible over }\C
 \quad\Longleftrightarrow\quad
 G\ \text{is $2$-transitive}.
\]
The polynomial $B_{P,Q}$ is symmetric, being the quotient of two
antisymmetric polynomials. Written in $u=z+w$ and $v=zw$, its irreducible
factors in $\C[u,v]$ correspond bijectively to the orbits of $G$ on unordered
pairs of distinct points. Thus the invariant Bezoutian is irreducible if and
only if $G$ is $2$-homogeneous. Individual components of $B_{P,Q}=0$ need not
be symmetric.
\end{proposition}

\begin{proof}
The ordered-pair assertion is the cited fiber-product theorem. It remains
to explain the passage to the invariant ring.

For the invariant statement, let $\sigma$ be the involution interchanging $z$
and $w$. If an irreducible factor $F$ of $B_{P,Q}$ satisfies
$F^\sigma=\lambda F$, then $\lambda=\pm1$, and $\lambda=-1$ is impossible,
because an antisymmetric polynomial vanishes on the diagonal and is divisible
by $z-w$. Hence a $\sigma$-stable irreducible factor can be chosen symmetric,
that is, it lies in $\C[z,w]^\sigma=\C[u,v]$, and it is clearly irreducible
there. If $F$ and $F^\sigma$ are distinct irreducible factors, their product
is symmetric and irreducible in $\C[u,v]$: by unique factorization in
$\C[z,w]$, any symmetric divisor of $FF^\sigma$ divisible by $F$ is also
divisible by $F^\sigma$. Thus the irreducible factors of $\overline B_{P,Q}$
correspond to the orbits of $\sigma$ on the set of nontrivial orbitals, that
is, to the orbits of $G$ on unordered pairs.
\end{proof}

Call an orbital \emph{connected} when its underlying undirected graph
on the generic fiber is connected. The standard component dictionary
then gives the following formulation for finite-orbit correspondences.

\begin{theorem}[Orbital dictionary]\label{th:dictionary}
Irreducible non-diagonal finite-orbit correspondences with generic orbit size
$N$ are in bijection with pairs $(f,\mathcal O)$, where $f$ is a rational
function of degree $N$, taken up to M\"obius transformations of the target,
and $\mathcal O$ is a connected nontrivial orbital of $\Mon(f)$. Under this
bijection, the graph of a generic orbit, with edges oriented from $z$ to $w$,
is the orbital digraph of $\mathcal O$. This digraph is the complete directed
graph, with both orientations of every edge, precisely when $\Mon(f)$ is
$2$-transitive; in that case the curve is the whole residual fiber curve.
\end{theorem}

\begin{proof}
The component associated with $(f,\mathcal O)$ has its orbits inside
fibers of $f$. Connectedness of the orbital graph makes a generic fiber
a single full orbit, so $f$ is its minimal quotient. Conversely, apply
\cref{th:intro-fio,prop:monodromy} to an irreducible finite-orbit
correspondence and its minimal quotient. The associated orbital must
be connected: otherwise its graph components form a nontrivial block
system and give a smaller right composition factor constant on the
orbits. These constructions are inverse. The complete-graph assertion
is the definition of $2$-transitivity.
\end{proof}

\begin{remark}[Relation with Poncelet curves]\label{rem:poncelet}
Mirman and Shukla characterize the Bezoutian equations associated with
their complex Poncelet curves \cite[Theorem~1.1]{MirmanShukla} and discuss
the quadratic case in \cite[\S4]{MirmanShukla}. Their completeness condition
concerns all pairs in a generic fiber. The whole residual fiber relation
has this condition whether or not it is irreducible. It is therefore
incorrect to identify their setting with $2$-transitive monodromy:
$2$-transitivity is the additional condition that this residual curve be
irreducible. Their geometric hypotheses must also be retained when passing
from a Bezoutian equation to a Poncelet realization.

\end{remark}

Equivalently, use a common Galois cover $Y$ with finite group $G$ and a
point stabilizer $H$, so $Y/H\simeq\Pone$; write $q:Y\to Y/H$ for
the quotient map. This is the orbital form of
the classical common-cover construction
\cite[Theorem~1 and condition~(H.3)]{FriedPoncelet}; compare
\cite[Lemma~4.2]{Krishnamoorthy} for hyperbolic correspondences.
Here the finite group is given, so hyperbolicity is unnecessary.

\begin{corollary}[Galois normal form]\label{cor:galois-normal}
A non-diagonal irreducible curve $\Gamma\subset\Pone\times\Pone$ is a
finite-orbit correspondence if and only if it has a parameterization
\[
 \Gamma=\overline{\{(q(y),q(\gamma y)):y\in Y\}}
\]
for data $(Y,G,H,\gamma)$ such that $G$ is finite and acts faithfully on the smooth projective
curve $Y$, the quotient $Y/H$ (and hence $Y/G$) has genus zero, $H$ is core-free in $G$,
$\gamma\notin H$, and $G$ is generated by $H$ together with the double coset
$H\gamma H$. The generic orbit size is then $[G:H]$, and both projection degrees of 
$\Gamma$ are equal to $|H\gamma H|/|H|$.
Interchanging $z$ and $w$ replaces $H\gamma H$ by $H\gamma^{-1}H$.
\end{corollary}

Here the generation condition expresses connectedness of the orbital
graph, and the double-coset index is its valency.

The dictionary also gives a simple necessary condition that can be read off
from an equation. Choose the target coordinate of the minimal quotient so that
$f(\infty)=\infty$, that is, $f=P/Q$ with $\deg P=N>\deg Q=\ell$. The top
homogeneous part of its residual fiber relation is then, up to a nonzero
scalar,
\[
 z^\ell w^\ell h_{N-\ell-1}(z,w),
 \qquad \text{where } h_j(z,w)=\frac{z^{j+1}-w^{j+1}}{z-w}.
\]

\begin{corollary}[Leading-form test]\label{cor:leading-test}
Let $F(z,w)=0$ be a reduced equation of a finite-orbit correspondence with no
diagonal component. Then the top homogeneous part of $F$ divides, up to a
nonzero scalar,
\[
 z^\ell w^\ell h_{N-\ell-1}(z,w)
 \quad\text{for some }N>\ell\ge0.
\]
In particular, every linear factor of the leading form not associated with $z$
or $w$ has the form $z-\zeta w$ with $\zeta$ a nontrivial $(N-\ell)$th root of
unity, and the leading form is squarefree away from $zw=0$.
\end{corollary}

\begin{proof}
By \cref{th:intro-fio}, $F$ divides the residual fiber relation of its minimal
quotient. Taking top homogeneous parts gives the assertion, and the linear
factors $z-\zeta w$ of $h_{N-\ell-1}$ are pairwise distinct.
\end{proof}

In projection degrees $(1,1)$, the irreducible finite-orbit correspondences are
exactly the graphs of finite-order M\"obius transformations; the graph
is symmetric precisely when the transformation is an involution or the
identity.

\section{Rational maps with a marked pole fiber}\label{sec:rational}

Let $V_d=H^0(\Pone,\mathcal O(d))$ be the space of binary forms of degree $d$.
A base-point-free pencil $\Lambda=\langle A,B\rangle\in\Gr(2,V_d)$ determines
the degree-$d$ map $A/B$ up to postcomposition with a M\"obius
transformation. The map
\begin{equation}\label{eq:pencil-embedding}
 \langle A,B\rangle\longmapsto[A(z)B(w)-A(w)B(z)]
\end{equation}
is the Pl\"ucker embedding followed by the natural inclusion
$\Lambda^2V_d\hookrightarrow V_d\otimes V_d$.
Dividing by the diagonal identifies this space with the symmetric forms
of bidegree $(d-1,d-1)$. Hence Pl\"ucker degrees are also degrees in the
projective space of symmetric Bezoutians.

We now fix a marked fiber. Every rational map of degree $k$ can be
written, after a M\"obius change of the target, as $f=P/Q$ with
$f(\infty)=\infty$ and coprime $P,Q$ of degrees $k>\ell$. Rescaling the pair
and scaling the target makes $P$ and $Q$ monic. The remaining target
translation replaces $P$ by $P+cQ$, so there is a unique choice with
$[t^\ell]P=0$. For $k>\ell\ge0$, let
$\mathcal P_{k,\ell}$ be the affine space of pairs $(P,Q)$ with $P,Q$ monic of
degrees $k,\ell$ and $[t^\ell]P=0$; it has dimension $k+\ell-1$. On the locus
where $P$ and $Q$ are coprime, the map $f=P/Q$ has the pole $\infty$ of order
\begin{equation}\label{eq:marked-index}
 e_f(\infty)=k-\ell,
\end{equation}
its finite poles are the roots of $Q$, and the top homogeneous part of
$B_{P,Q}$ is $z^\ell w^\ell h_{k-\ell-1}(z,w)$, with $h_j$ as in
\cref{cor:leading-test}.

The factorization under composition is classical; compare
\cite[\S3.2]{Pak2011}. We include the homogeneous identity and the
proper parameter space, which will also be used in the degree calculation.

\begin{proposition}[Composition and closedness]\label{prop:rational-composition}
Let $f=g\circ h$ with $\deg g=r>1$ and $\deg h=s>1$, and write $h=A/B$ and
$g=C/D$ in coprime form. Let $\widehat B_g$ denote the bihomogenization of
$B_g=B_{C,D}$ to bidegree $(r-1,r-1)$. Then, up to a nonzero scalar,
\begin{equation}\label{eq:rational-composition}
 B_f(z,w)=B_h(z,w)\,
 \widehat B_g\bigl((A(z),B(z)),(A(w),B(w))\bigr).
\end{equation}
Both factors are nonconstant and symmetric. On the space of base-point-free
pencils of degree $rs$, the condition of having a right composition factor of
degree $s$ is closed, and a given map has only finitely many right composition
factors up to postcomposition.
\end{proposition}
\begin{proof}
For $X=(X_1,X_2)$, let $\widehat C(X)=X_2^r C(X_1/X_2)$ and
$\widehat D(X)=X_2^r D(X_1/X_2)$ be the degree-$r$ homogenizations of $C$ and
$D$. By bihomogenizing the defining identity for $B_g$, we have
\[
 \widehat C(X)\widehat D(Y)-\widehat C(Y)\widehat D(X)
 =(X_1Y_2-X_2Y_1)\widehat B_g(X,Y)
\]
for $Y=(Y_1,Y_2)$. Substituting $(A(z),B(z))$ and $(A(w),B(w))$ and dividing
by $z-w$ gives \eqref{eq:rational-composition} exactly for the coprime
representative $(\widehat C(A,B),\widehat D(A,B))$ of $f$. The factors are
symmetric. Their bidegrees are bounded by $(s-1,s-1)$ and $(rs-s,rs-s)$,
respectively. The coefficient of $z^{rs}$ in $P(z)Q(w)-P(w)Q(z)$ is
nonzero for a coprime representative $(P,Q)$ of degree $rs$, so $B_f$ has
bidegree $(rs-1,rs-1)$ and both bounds are attained.

In affine coordinates the second factor is
\[
 B(z)^{r-1}B(w)^{r-1}B_g(h(z),h(w)).
\]
For polynomial $h$ represented with $B=1$, this is simply $B_g(h(z),h(w))$.

For closedness, let $\cS$ be the tautological subbundle on $\Gr(2,V_s)$.
Multiplication of forms embeds $\Sym^r\cS$ into the trivial bundle with fiber
$V_{rs}$: for linearly independent $A,B$, the ratio $A/B$ is nonconstant, so
the products $B^r,AB^{r-1},\ldots,A^r$ are linearly independent. This
embedding induces a morphism $\Gr_2(\Sym^r\cS)\to\Gr(2,V_{rs})$, which is
proper because $\Gr_2(\Sym^r\cS)$ is. Its image is therefore closed. A point
of $\Gr_2(\Sym^r\cS)$ is a pair consisting of an inner pencil
$\langle A,B\rangle$ and an outer pencil
$\langle\widehat C(A,B),\widehat D(A,B)\rangle$. If the image pencil is
base-point free, then so is $\langle A,B\rangle$, since otherwise every form in
the image vanishes on the common base divisor raised to the $r$th power; and
so is $\langle\widehat C,\widehat D\rangle$, since a common zero $(a:b)$ of
$\widehat C$ and $\widehat D$ would give a common zero of the image forms at
every point where $(A:B)=(a:b)$, and such points exist because $h$ is
surjective. Hence the intersection of the image with the base-point-free locus
is exactly the locus of maps with a right composition factor of degree $s$.
Finally, right composition factors up to postcomposition correspond to block
systems of the finite group $\Mon(f)$, so there are finitely many of them.
\end{proof}

We record two elementary permutation-group facts. The block calculation
is standard. The second lemma is a relative form of the familiar argument
using an isolated simple branch point to obtain a full symmetric kernel.
Related inertia arguments for iterates and compositions occur in
\cite{Odoni,Juul,JuulKuMaTu,KonigNeRo}.

\begin{lemma}[Wreath products]\label{lem:blocks}
The imprimitive action of $S_s\wr S_r$ on $rs$ points, with $r,s\ge2$, has
exactly one nontrivial block system, namely the given one, with $r$ blocks of
size $s$.
\end{lemma}
This is the usual block calculation for the imprimitive wreath product:
the point stabilizer has orbits of sizes $1$, $s-1$, and $s(r-1)$;
only $1$, $s$, and $rs$ among their unions containing the point can be
block sizes. See \cite{Wielandt} for the block-system correspondence.

\begin{lemma}[Wreath criterion]\label{lem:wreath}
Let $f=g\circ h$ with $\deg g=r>1$ and $\deg h=s>1$, and assume that $h$ is
indecomposable. If some critical value of $f$ has exactly one ramified
preimage, and this preimage has local degree two, then
$\Mon(f)=S_s\wr\Mon(g)$ in its imprimitive action. In particular, if
$\Mon(g)=S_r$, then $\Mon(f)=S_s\wr S_r$.
\end{lemma}
\begin{proof}

Let $G=\Mon(f)$. The fibers of $h$ over the points of a generic fiber of $g$
form a system of $r$ blocks $\Omega_1,\dots,\Omega_r$ of size $s$, so $G$ is a
subgroup of $\Sym(\Omega_1)\wr S_r$, and the group induced on the set of blocks
is $\Mon(g)$. The setwise stabilizer of a block $\Omega_i$ acts on $\Omega_i$
through $\Mon(h)$, which is primitive because $h$ is indecomposable.

By hypothesis, the branch cycle of $f$ over the given critical value is a
transposition $\tau$. It preserves every block, because a permutation in the
wreath product that moves some block moves at least $2s>2$ points. Consider the
graph on the fiber whose edges are the transpositions conjugate to $\tau$ in
$G$. All its edges lie inside blocks, and its connected components form a
$G$-invariant partition. Restricted to a block $\Omega_i$, this partition is
invariant under the block stabilizer, hence trivial or equal to $\{\Omega_i\}$
by primitivity. Since $\tau$ is an edge inside some block and $G$ permutes the
blocks transitively, every block contains an edge, so the graph is connected on
every block. Transpositions forming a connected graph on a set generate the full
symmetric group of that set. Therefore the conjugates of $\tau$ generate the
subgroup
\[
 \Sym(\Omega_1)\times\cdots\times\Sym(\Omega_r)\simeq S_s^r,
\]
which is the kernel of the action of $S_s\wr S_r$ on the blocks. Since $G$
contains this kernel and maps onto $\Mon(g)$, we get $G=S_s\wr\Mon(g)$.
\end{proof}

\subsection{Simple finite poles}

\begin{theorem}[Dichotomy for simple finite poles]\label{th:separable-rational}
Let $(P,Q)\in\mathcal P_{k,\ell}$ with $1\le\ell\le k-2$, where $P,Q$ are
coprime and $Q$ is squarefree. Then either $f=P/Q$ is decomposable, or
$\Mon(f)$ is $2$-transitive. Consequently, on this locus the following are
equivalent: $f$ is decomposable; $B_f$ is reducible; $\overline B_f$ is
reducible.
\end{theorem}
\begin{proof}

Suppose that $f$ is indecomposable, so that $\Mon(f)$ is primitive.
The marked pole fiber supplies a $(k-\ell)$-cycle fixing $\ell$ points.
Jordan's theorem, in the form recalled by Jones \cite[\S2, p.~161]{Jones},
says that a primitive group containing such a cycle is
$(\ell+1)$-transitive. Since $\ell\ge1$, the group is $2$-transitive.
This applies uniformly for all values of $\ell$ in the statement; the
classification of primitive groups containing cycles is not required.

If $f$ is decomposable, \eqref{eq:rational-composition} shows that
$\overline B_f$ is reducible, and this implies that $B_f$ is reducible. If
$B_f$ is reducible, then $\Mon(f)$ is not $2$-transitive by
\cref{prop:monodromy}, so $f$ is decomposable by the first part.
\end{proof}

For the remainder of this subsection assume $1\le\ell\le k-2$, and set
\[
 \mathcal U_{k,\ell}=\{(P,Q)\in\mathcal P_{k,\ell}:\ \gcd(P,Q)=1,
 \ Q\text{ squarefree}\}.
\]
For $s\mid k$ with $1<s<k$, let $\mathcal D_s\subset\mathcal U_{k,\ell}$ be the
locus of maps having a right composition factor of degree $s$.

\begin{theorem}[Composition components]\label{cor:rational-strata}
Let $r=k/s$ and $b=k-\ell$. The locus $\mathcal D_s$ is nonempty precisely in
the following two cases:
\begin{center}
\begin{tabular}{@{}clll@{}}
\toprule
case&condition&dimension&codimension\\
\midrule
{\rm (i)}&$s\ge b$&$2r+2s-b-3$&$2(r-1)(s-1)$\\
{\rm (ii)}&$s\mid\ell$, $a:=b/s$&$2r+s-a-2$&$k+\ell+1-2r-s+a$\\
\bottomrule
\end{tabular}
\end{center}
The two cases overlap only for $s=b$, where they describe the same family.
Every nonempty $\mathcal D_s$ is closed in $\mathcal U_{k,\ell}$, irreducible
and rational. Its general member has monodromy $S_s\wr S_r$. Its invariant
Bezoutian has exactly two irreducible factors, which have bidegrees
$(s-1,s-1)$ and $(k-s,k-s)$ after substituting $u=z+w$ and $v=zw$.
The nonempty loci $\mathcal D_s$ are exactly the irreducible
components of the locus in $\mathcal U_{k,\ell}$ where $\overline B_f$ is
reducible.
\end{theorem}
\begin{proof}
Let $f=g\circ h$ with $\deg h=s$. After an intermediate M\"obius
transformation we may assume $h(\infty)=\infty$. Then
$e_g(\infty)e_h(\infty)=e_f(\infty)=b$, while every other point of the pole
fiber of $f$ is a simple pole. If $e_g(\infty)>1$, then every point of $h^{-1}(\infty)$ 
is a pole of $f$ of order at least $e_g(\infty)$, so $h^{-1}(\infty)=\{\infty\}$. 
Hence either $e_g(\infty)=1$, which is Case (i), or $h^{-1}(\infty)=\{\infty\}$, 
which is Case (ii); the two overlap when $h$ is a polynomial and $e_g(\infty)=1$, 
that is, when $s=b$. In Case (i), $h=A/B$ with numerator and denominator of degrees $s$
and $s-b$, and $g=C/D$ with degrees $r$ and $r-1$. In Case (ii), $h=A$ is a
polynomial of degree $s$, $e_g(\infty)=a=b/s$, and $g=C/D$ has degrees $r$ and
$r-a$, so that $\ell=(r-a)s$. These conditions are also sufficient: one chooses
the finite poles of $g$ simple and away from the critical values of $h$, and,
in Case (i), the polynomial $B$ squarefree. If both conditions hold, then
$s\mid b$ and $s\ge b$, hence $s=b$, and both cases describe the polynomial
right composition factors of degree $b$.

We normalize $A$ and $B$ to be monic with $[t^{\deg B}]A=0$, as for the marked
pole normalization above. This fixes the intermediate coordinate, and in
Case (ii) it means $A(0)=0$. Write $g(y)=C(y)/D(y)$, where $y=h(t)$ is the
intermediate coordinate, and take $C$ and $D$ monic. If $m=\deg D$ and
$\delta=\deg B$, then $\ell=ms+(r-m)\delta$. Writing $c_j=[y^j]C$, we have
$P=\sum_{j=0}^r c_j A^jB^{r-j}$. The monic polynomial $A^jB^{r-j}$ has
degree $js+(r-j)\delta$, strictly increasing in $j$. Thus $[t^\ell]P=0$
determines $[y^m]C$, since $[t^\ell](A^mB^{r-m})=1$. Here $m=r-1$ in
Case (i) and $m=r-a$ in Case (ii). The remaining parameters range over
a nonempty open subset of an affine space of dimension
\[
 (2s-b-1)+(2r-2)\quad\text{in Case (i)},\quad (s-1)+(2r-a-1)\quad\text{in Case (ii)}.
\]
The fibers of the parameter map onto $\mathcal D_s$ are finite, since the
normalization fixes a chosen right composition factor uniquely and there are
only finitely many such factors. This proves irreducibility and the dimension
formulas; the codimensions follow from $\dim\mathcal P_{k,\ell}=k+\ell-1$.
Closedness follows from \cref{prop:rational-composition}.

For $2\le b'\le d$, we construct a degree-$d$ cover with one branch cycle
of type $(b',1^{d-b'})$ and simple branching over every other branch value.
With permutations composed from right to left, the inverse of
$(1\,2\,\cdots\,b')$ is $(1\,2)(1\,3)\cdots(1\,b')$. List the cycle, then
these $b'-1$ transpositions, then two consecutive copies of $(1\,j)$ for each
$b'<j\le d$. Each pair represents two distinct branch values and has
product one. The tuple has product one, generates $S_d$ through
the transpositions $(1\,j)$ for $2\le j\le d$, and has index sum
$(b'-1)+(b'-1)+2(d-b')=2d-2$. Riemann's existence theorem and Riemann--Hurwitz
therefore give a connected cover of genus zero. For an unramified marked
fiber, use the tuple consisting only of the pairs $(1\,j),(1\,j)$ for
$2\le j\le d$, with the same product, generation, and index properties.

First choose $g$ by this construction with its prescribed pole profile.
Then choose $h$ with its prescribed pole profile and with distinct simple
branch values away from $\infty$ and $g^{-1}(\{\infty\}\cup\operatorname{Br}(g))$,
where $\operatorname{Br}(g)$ is the set of branch values of $g$. We can also
choose these values to have distinct images under $g$, since each choice
excludes only finitely many further points. There is at least one such value
for $h$, since the marked fiber contributes at most $s-1$ to the total index
$2s-2$. Over its image under $g$, the composite has exactly one ramified
point of local degree two. The simple branching and separation conditions
define a nonempty Zariski-open subset of the
irreducible parameter space. There $\Mon(h)=S_s$ and $\Mon(g)=S_r$, since
the simple branch cycles generate the monodromy of each map (the marked
branch cycle is their inverse product when the marked fiber is ramified).
\Cref{lem:wreath} now gives $\Mon(g\circ h)=S_s\wr S_r$ on this dense open
subset.

The wreath product $S_s\wr S_r$ has exactly two orbits on unordered pairs of
distinct points: pairs inside one block and pairs in different blocks. By
\cref{prop:monodromy} this gives the two invariant factors, of the stated
bidegrees. By \cref{lem:blocks}, the generic member of $\mathcal D_s$ has a
unique right composition factor of degree $s$, so the parameter map is
birational and $\mathcal D_s$ is rational. Moreover, a member with full wreath
monodromy lies in no $\mathcal D_{s'}$ with $s'\ne s$, so none of these closed
irreducible loci contains another. By \cref{th:separable-rational}, their union
is the reducibility locus.
\end{proof}

\begin{example}\label{ex:rational-six-two}
For $(k,\ell)=(6,2)$ only Case (ii) with $s=2$ occurs. The locus
$\mathcal D_2$ has dimension $4$ and codimension $3$. An explicit member is
\[
 f(t)=\frac{t^6}{t^2-c}=\left(\frac{y^3}{y-c}\right)\circ t^2,
 \qquad c\ne0,
\]
with
\[
 B_f=(z+w)\bigl(z^4w^2+z^2w^4-c(z^4+z^2w^2+w^4)\bigr).
\]
For $(k,\ell)=(6,1)$ neither case occurs, so every coprime pair in
$\mathcal P_{6,1}$ has an absolutely irreducible Bezoutian.
\end{example}

\begin{remark}\label{rem:off-separable}
Both hypotheses of \cref{th:separable-rational} are needed.
The decomposition locus remains closed without the assumption that $Q$ is
squarefree, but then the ramification patterns of
\cref{cor:rational-strata} no longer exhaust it; moreover, indecomposable maps
with reducible Bezoutian appear. For instance, if $f$ is a Dickson map or an
elliptic quotient map of prime degree $p\ge5$ (see \cref{sec:prime}) and the
marked fiber is chosen to be a fiber of type $(2,\ldots,2,1)$ with $\infty$ at
a ramified point, then $\ell=p-2$ and $Q$ has double roots.

The range $\ell\le k-2$ is also essential. For $\ell=k-1$ the point $\infty$
is unramified, and the argument above gives no information. For example,
\[
 P(t)=t^5,\qquad Q(t)=\frac{t^5-(t-1)^5}{5}=t^4-2t^3+2t^2-t+\tfrac15
\]
is a coprime pair in $\mathcal P_{5,4}$, and $Q$ is squarefree, its roots being
$\zeta/(\zeta-1)$ with $\zeta^5=1\ne\zeta$. With the M\"obius transformation
$h(t)=t/(t-1)$ one has
\[
 \frac{P}{Q}=\frac{5h^5}{h^5-1},
\]
so $P/Q$ is M\"obius equivalent to $t^5$. Its monodromy is cyclic of order
five, $B_{P,Q}$ splits into four factors of bidegree $(1,1)$, and
$\overline B_{P,Q}$ is the product of two irreducible factors, corresponding to
the two pairs $\{\zeta,\zeta^{-1}\}$; yet $P/Q$ is indecomposable, its degree
being prime. We do not attempt a complete classification outside the range of
\cref{th:separable-rational}.
\end{remark}

\subsection{Absolute and invariant factors}

\begin{corollary}[Kantor's obstruction]\label{cor:kantor}
If $B_f$ is reducible and $\overline B_f$ is irreducible, then
$\deg f=q$ is a prime power with $q\equiv3\pmod4$, and $\Mon(f)$ is a
$2$-homogeneous subgroup of $\operatorname{A\Gamma L}_1(\mathbb F_q)$.
Moreover,
\begin{equation}\label{eq:paired-factors}
 B_f(z,w)=c\,\Phi(z,w)\Phi(w,z),
 \qquad \deg_z\Phi=\deg_w\Phi=(q-1)/2,
\end{equation}
where $c\in\C^*$ and the two factors $\Phi(z,w)$ and $\Phi(w,z)$ are distinct
and irreducible.
\end{corollary}
\begin{proof}

By \cref{prop:monodromy}, $\Mon(f)$ is $2$-homogeneous but not
$2$-transitive. In degree three the assertion is immediate; in larger
degree it follows from Kantor's classification
\cite[Theorem~1, p.~261]{Kantor}. Such a group has exactly two nontrivial orbitals,
interchanged by transposition, each of size $q(q-1)/2$; by
\cref{prop:monodromy} they give the two factors.
\end{proof}

\begin{remark}\label{rem:kantor-semilinear}
Kantor's theorem gives containment in a semilinear affine group; it does
not say that the scalar subgroup contains all square multipliers. In
\cref{th:all-discrepancy} we impose the genus-zero condition without making
this extra assumption and obtain a classification in all degrees.

\end{remark}

\section{Invariant factorization for polynomial maps}\label{sec:polynomial}

Let $P\in\C[t]$ have degree $k\ge2$. Since $B_P(z,w)$ is symmetric, there is a
unique polynomial $\overline B_P\in\C[u,v]$ with
$B_P(z,w)=\overline B_P(z+w,zw)$. We assign the weights $\wt(u)=1$ and
$\wt(v)=2$, so that $\overline B_P$ has weighted degree $k-1$. The
\emph{invariant factorization type} of $P$ is the multiset of weighted degrees
of the irreducible factors of $\overline B_P$ in $\C[u,v]$; we display it as a
tuple whose ordering is immaterial. Since an orbit on unordered pairs may
combine two orbitals that are interchanged by transposition, invariant and
absolute factorization need not coincide.

\subsection{The complete criterion}

We normalize the Dickson polynomials $D_n(t,a)$ by
$D_n(x+a/x,a)=x^n+(a/x)^n$, so that
\[
 D_n(t,a)=\sum_{i=0}^{\lfloor n/2\rfloor}\frac{n}{n-i}\binom{n-i}{i}(-a)^it^{n-2i},
 \qquad D_n(t,0)=t^n,
\]
and $D_{rs}(t,a)=D_r(D_s(t,a),a^s)$. We use the following elementary
observation and a classical identity.

\begin{lemma}\label{lem:linear-in-v}
Let $g(u)\in\C[u]$ and $\lambda\in\C^*$. Then $g(u)-\lambda v$ is
irreducible in $\C[u,v]$.
\end{lemma}

\begin{proof}
Substitution $v=\lambda^{-1}g(u)$ identifies
$\C[u,v]/(g(u)-\lambda v)$ with the integral domain $\C[u]$.
Thus $g(u)-\lambda v$ generates a nonzero prime ideal and is irreducible.
\end{proof}

The second fact is the classical factorization of Dickson difference
quotients; see \cite[Proposition~1.7, p.~318]{Turnwald}.

\begin{lemma}[Dickson difference identity]\label{lem:dickson-difference}
Let $n\ge2$, let $\zeta$ be a primitive $n$th root of unity, and put
$\alpha_j=\zeta^j+\zeta^{-j}$ and $\beta_j=\zeta^j-\zeta^{-j}$. Then
\begin{equation}\label{eq:dickson-odd}
 \frac{D_n(z,a)-D_n(w,a)}{z-w}
 =(z+w)^{\varepsilon}\prod_{j=1}^{\lfloor (n-1)/2\rfloor}
 \bigl(z^2-\alpha_jzw+w^2+\beta_j^2a\bigr),
\end{equation}
where $\varepsilon=0$ if $n$ is odd and $\varepsilon=1$ if $n$ is even.
\end{lemma}

\begin{theorem}[Invariant polynomial criterion]\label{th:fried}
Let $P\in\C[t]$ have degree $k\ge3$. The invariant Bezoutian $\overline B_P$
is reducible in $\C[u,v]$ if and only if either $P$ is decomposable, or
$k=p\ge5$ is prime and $P$ is affinely equivalent to $D_p(t,a)$ for some
$a\in\C$.
\end{theorem}

\begin{proof}
If $P=G\circ H$ is decomposable, then
\begin{equation}\label{eq:composition-bez}
 B_P(z,w)=B_H(z,w)B_G(H(z),H(w))
\end{equation}
is a product of nonconstant symmetric polynomials. Affine changes of source
and target preserve both absolute and invariant reducibility. In fact,
an affine change of source acts by the same invertible substitution on
$z$ and $w$, which commutes with their interchange, and the changes affect
the Bezoutian by a nonzero scalar in addition to the substitution.
Invariant reducibility implies absolute reducibility. For indecomposable
$P$, Fried's theorem
\cite[Theorem~1, p.~45]{FriedSchur}, in the precise form of
\cite[Theorem~1, p.~313]{Turnwald}, says that absolute reducibility
forces $P$ to be affinely equivalent to a prime-degree Dickson polynomial,
with the cyclic case the only possibility in degree three.
A normalized cubic has invariant Bezoutian $u^2-v+a$, which is irreducible.
For prime $p\ge5$, \eqref{eq:dickson-odd} gives $(p-1)/2$ distinct
invariant factors
$u^2-(\alpha_j+2)v+\beta_j^2a$. Each is irreducible because
$\alpha_j+2\ne0$, so its coefficient of $v$ is a nonzero constant.
\end{proof}

In prime degree $k=p\ge5$, consider the representative $P=D_p(t,a)$.
For $a=0$ we have $P=D_p(t,0)=t^p$, and each invariant quadratic is the product
of the linear factors $z-\zeta^jw$ and $z-\zeta^{-j}w$, whose orbitals are
interchanged by transposition. For $a\ne0$, each conic
$z^2-\alpha_jzw+w^2+\beta_j^2a=0$ is nondegenerate, hence already irreducible
in $\C[z,w]$, and the corresponding orbitals are self-paired.

\subsection{The normalized slice and the composition loci}

Affine changes of the variable $t$ act on $\C[u,v]$ by automorphisms, and
affine changes of the target change $B_P$ only by a nonzero scalar. It is
therefore natural to work in the normalized slice
\[
 \widetilde{\mathcal P}_k
 =\{t^k+a_{k-2}t^{k-2}+\cdots+a_1t\}\simeq\A^{k-2}
\]
of monic, centered polynomials with zero constant term. This slice removes
the translations on both sides, but not all of affine equivalence: a residual
action of $\Gm$ remains,
\begin{equation}\label{eq:residual-scaling}
 P(t)\longmapsto \alpha^{-k}P(\alpha t),\qquad\alpha\in\C^* ,
\end{equation}
and every locus considered below is invariant under it. On the Dickson family
it acts by $a\mapsto a\alpha^{-2}$. For $k=rs$ with $r,s>1$, let $\cD_{r,s}$ be
the image of the map $(G,H)\mapsto G\circ H$ defined on pairs of polynomials
of the form
\begin{equation}\label{eq:normalized-composition}
 H(t)=t^s+\sum_{j=1}^{s-2}c_jt^j,
 \qquad
 G(t)=t^r+b_{r-1}t^{r-1}+\cdots+b_1t.
\end{equation}

The next statement is the centered form of the established geometry of
decomposition varieties. The polynomial recovery map is the
Newton--Taylor construction; see \cite{vonzurGathen1990} and
\cite[Theorem~2.2]{vzGMat2017}. We recall it to fix the normalization.

\begin{proposition}[Composition loci]\label{prop:composition-loci}
The normalized composition map
\[
 \iota_{r,s}\colon\ \widetilde{\mathcal G}_r\times\widetilde{\mathcal H}_s
 \simeq\A^{r-1}\times \A^{s-2}\ \longrightarrow\ \widetilde{\mathcal P}_k,
\]
where $\widetilde{\mathcal G}_r$ and $\widetilde{\mathcal H}_s$ are the spaces
of polynomials $G$ and $H$ as in \eqref{eq:normalized-composition}, is a closed
embedding. Consequently,
\[
 \cD_{r,s}\simeq\A^{r+s-3},\qquad
 \operatorname{codim}\cD_{r,s}=(r-1)(s-1),
\]
and every point of $\cD_{r,s}$ has a unique decomposition of the form
\eqref{eq:normalized-composition}.
\end{proposition}

The inverse is given by the polynomial Newton--Taylor recovery map in the
cited references. The centered and zero-constant conditions simply remove
the two translation parameters from their monic normalization.

The standard generic wreath-product monodromy gives the following
unordered-pair consequence.

\begin{theorem}[Generic factorization]\label{th:generic-composition}
There is a dense open subset $\cD_{r,s}^{\circ}\subset\cD_{r,s}$ on which
$\Mon(P)=S_s\wr S_r$, acting on $r$ blocks of size $s$. For
$P=G\circ H\in\cD_{r,s}^{\circ}$, the two factors in
\eqref{eq:composition-bez} are irreducible in $\C[u,v]$, of weighted degrees
$s-1$ and $rs-s$.
\end{theorem}

\begin{proof}
Choose $G,H$ with simple critical points and separated critical values,
including separation of the images under $G$ of the critical values of
$H$. Then $\Mon(G)=S_r$ and $\Mon(H)=S_s$, and an isolated simple branch
point supplies the transposition required by \cref{lem:wreath}.
The two unordered-pair orbits of $S_s\wr S_r$ consist of pairs in the
same block and pairs in distinct blocks. They give the two stated factors
and their degrees by \cref{prop:monodromy}.
\end{proof}

The decomposition components are known from
\cite{vonzurGathen2015,vzGMat2017}. Combining them with Fried's criterion
and the invariant Dickson factorization gives the following description.

\begin{corollary}[The normalized reducibility locus]\label{cor:complete-locus}
Let $\widetilde{\mathcal R}_k\subset\widetilde{\mathcal P}_k$ be the locus
where $\overline B_P$ is reducible. Then
\begin{equation}\label{eq:complete-locus}
 \widetilde{\mathcal R}_k
 =\bigcup_{\substack{k=rs\\r,s>1}}\cD_{r,s}\ \cup\ \mathcal X_k,
\end{equation}
where $\mathcal X_k=\{D_k(t,a)-D_k(0,a):a\in\C\}$ if $k\ge5$ is prime, and
$\mathcal X_k=\emptyset$ otherwise. The sets $\cD_{r,s}$ and $\mathcal X_k$
appearing in \eqref{eq:complete-locus} are precisely the irreducible components
of $\widetilde{\mathcal R}_k$. In composite degree, the normalized Dickson
family $\{D_k(t,a)-D_k(0,a)\}$ is contained in every $\cD_{r,s}$.
\end{corollary}

\begin{proof}
The known normalization of decompositions and \cref{th:fried} give the
set-theoretic equality. A generic point of $\cD_{r,s}$ has monodromy
$S_s\wr S_r$, whose only nontrivial block system has block size $s$.
Thus no component is contained in another. In prime degree, centered
affine copies of the Dickson polynomial give precisely $\mathcal X_k$;
its coefficient of $t^{k-2}$ is $-ka$, so this is an affine line.
For composite $k=rs$, writing $c=D_s(0,a)$, the normalized factors
\[
 H=D_s(t,a)-c,\qquad G=D_r(t+c,a^s)-D_r(c,a^s)
\]
give $G\circ H=D_k(t,a)-D_k(0,a)$, proving the last assertion.
\end{proof}

If $p$ is the smallest prime divisor of a composite number $k$, the components
of largest dimension correspond to the most unbalanced factorizations, and
\[
 \dim\widetilde{\mathcal R}_k=k/p+p-3,
 \qquad
 \operatorname{codim}\widetilde{\mathcal R}_k=(p-1)(k/p-1).
\]

\subsection{The classical Dickson factorization in invariant coordinates}
The following proposition is a reformulation of the classical difference
identity, not an additional factorization theorem.

\begin{proposition}\label{prop:dickson-split}
For $k\ge3$ and $a\in\C$, put $P_{k,a}=D_k(t,a)-D_k(0,a)$. Its invariant
factorization type is $(2,\ldots,2)$ for $k$ odd and $(1,2,\ldots,2)$ for $k$
even. More precisely, with $\zeta$ a primitive $k$th root of unity and
$\alpha_j,\beta_j$ as in \cref{lem:dickson-difference},
\[
 \overline B_{P_{k,a}}(u,v)=
 \begin{cases}
 \displaystyle\prod_{j=1}^{(k-1)/2}
  \bigl(u^2-(\alpha_j+2)v+\beta_j^2a\bigr),&k\text{ odd},\\[2mm]
 \displaystyle u\prod_{j=1}^{(k-2)/2}
  \bigl(u^2-(\alpha_j+2)v+\beta_j^2a\bigr),&k\text{ even}.
 \end{cases}
\]
\end{proposition}

Indeed, substitute $z^2+w^2=u^2-2v$ in
\eqref{eq:dickson-odd}. Each displayed quadratic is irreducible in
$\C[u,v]$ because its coefficient of $v$ is nonzero and constant.

\section{Degrees of polynomial and rational composition loci}\label{sec:scroll}

The components of $\widetilde{\mathcal R}_k$ are affine spaces, and it is
natural to ask for the degrees of their projective closures. Throughout this
section, $k=rs$ with $r,s\ge2$. Let
\[
 W=\langle\, t,\ t^2,\ \ldots,\ t^{k-2},\ t^{k}\,\rangle\subset\C[t],
 \qquad \dim W=k-1,
\]
and identify $\mathbb P^{k-2}=\mathbb P(W)$; then $\widetilde{\mathcal P}_k$ is
the affine chart in which the coefficient of $t^k$ equals one. We write
$\overline{\cD}_{r,s}$ and $\overline{\mathcal X}_p$ for the closures in
$\mathbb P^{k-2}$ of the composition loci and of the Dickson curve.

\subsection{The scroll}

Give $\mathbb P^{s-2}$ the homogeneous coordinates $(c_0:c_1:\cdots:c_{s-2})$
and put
\[
 \widehat H_c(t)=c_0t^s+c_{s-2}t^{s-2}+\cdots+c_1t ,
\]
the homogenization of the right composition factor in
\eqref{eq:normalized-composition}.

\begin{lemma}\label{lem:scroll-bundle}
For every $c\in\mathbb P^{s-2}$ and every $1\le j\le r$, the polynomial
$\widehat H_c^{\,j}$ lies in $W$, and
$\widehat H_c,\widehat H_c^{\,2},\ldots,\widehat H_c^{\,r}$ are linearly
independent. Consequently, the map
\[
 \cE:=\mathcal O(-1)\oplus\cdots\oplus\mathcal O(-r)
 \xrightarrow{\ \varepsilon\ } W\otimes\mathcal O_{\mathbb P^{s-2}},
\]
sending the $j$th summand to the line spanned by $\widehat H_c^{\,j}$, is a
fiberwise injective map of vector bundles, and $\cE$ is a subbundle of rank $r$
of the trivial bundle $W\otimes\mathcal O_{\mathbb P^{s-2}}$.
\end{lemma}

\begin{proof}
The powers of $\widehat H_c$ have zero constant term. For $j<r$ their degrees
are at most $k-s\le k-2$, while $\widehat H_c^{\,r}$ has no term in $t^{k-1}$
because $\widehat H_c$ has no term in $t^{s-1}$. Thus all these powers lie in
$W$. Since $\widehat H_c$ is nonconstant, the powers have distinct degrees and
are linearly independent. Their coefficients are homogeneous of degree $j$ in
$c$, which gives the twists by $\mathcal O(-j)$ and a fiberwise injective
bundle map.
\end{proof}

\begin{theorem}[Polynomial composition degree]\label{th:scroll}
The morphism
\[
 \Psi\colon\ \mathbb P(\cE)\longrightarrow\mathbb P(W)=\mathbb P^{k-2}
\]
induced by $\varepsilon$ is birational onto its image, which is
$\overline{\cD}_{r,s}$. Thus $\overline{\cD}_{r,s}$ is swept out by a family of
$(r-1)$-dimensional linear spaces parameterized by $\mathbb P^{s-2}$, and
\begin{align*}
 \deg\overline{\cD}_{r,s}
 &=\int_{\mathbb P^{s-2}}s_{s-2}(\cE)
 =\bigl[h^{s-2}\bigr]\prod_{j=1}^{r}\frac{1}{1-jh}\\
 &=h_{s-2}(1,2,\ldots,r)=\Stir{r+s-2}{r},
\end{align*}
where $h_{m}$ denotes the complete homogeneous symmetric polynomial of degree
$m$ and $\Stir{\cdot}{\cdot}$ is the Stirling number of the second kind.
\end{theorem}

\begin{proof}
Over the open set $c_0=1$, write $H=\widehat H_c$. The fiber of $\mathbb P(\cE)$
maps onto $\mathbb P\langle H,H^2,\ldots,H^r\rangle$, whose intersection with the
chart $\widetilde{\mathcal P}_k$ of monic polynomials of degree $k$ is
\[
 \Bigl\{H^r+\sum_{j=1}^{r-1}b_jH^j\Bigr\}.
\]
Conversely, a polynomial of degree $k=rs$ in the image must come from a point
with $c_0\ne0$ and have a nonzero coefficient of $H^r$. Hence the intersection
of the image with $\widetilde{\mathcal P}_k$ is exactly $\cD_{r,s}$. Since
$\mathbb P(\cE)$ is irreducible and $\Psi$ is proper, the image is the closure
$\overline\cD_{r,s}$, and the uniqueness of normalized decompositions
(\cref{prop:composition-loci}) shows that $\Psi$ is birational onto it.

By \cref{lem:scroll-bundle},
$\mathcal O_{\mathbb P(\cE)}(-1)\subset\pi^*\cE\subset
W\otimes\mathcal O_{\mathbb P(\cE)}$, where $\pi\colon\mathbb P(\cE)\to\mathbb P^{s-2}$
is the projection, so
$\Psi^*\mathcal O_{\mathbb P^{k-2}}(1)=\mathcal O_{\mathbb P(\cE)}(1)$; let
$\xi$ be its first Chern class. Since $\dim\mathbb P(\cE)=r+s-3$, birationality
and the projection formula give
\[
 \deg\overline{\cD}_{r,s}
 =\int_{\mathbb P(\cE)}\xi^{\,r+s-3}
 =\int_{\mathbb P^{s-2}}\pi_*\bigl(\xi^{(r-1)+(s-2)}\bigr)
 =\int_{\mathbb P^{s-2}}s_{s-2}(\cE),
\]
where $s(\cE)=c(\cE)^{-1}$ is the Segre class \cite[Ch.~3]{Fulton}. Since
$c(\cE)=\prod_{j=1}^r(1-jh)$, with $h$ the hyperplane class of
$\mathbb P^{s-2}$, we get
$s_{s-2}(\cE)=[h^{s-2}]\prod_{j=1}^{r}(1-jh)^{-1}=h_{s-2}(1,\ldots,r)$.
Finally, the standard generating function
$\sum_{m\ge0}\Stir{m+r}{r}x^m=\prod_{j=1}^r(1-jx)^{-1}$ for Stirling numbers of
the second kind \cite[Ch.~1]{Stanley} identifies this coefficient with
$\Stir{r+s-2}{r}$.
\end{proof}

\subsection{Uncentered polynomials and interpolation}

The choice of normalization affects the projective degree. Let $d,e\ge1$, $n=de$, and
let $C_{n,d}\subset\A^n$ denote the monic degree-$n$ polynomials
$F=G\circ H$, where $G$, $H$ are monic of degrees $d,e$ and $H(0)=0$.
The constant coefficient of $G$ is free. This is the normalization of
\cite[\S3]{vzGMat2024}.

\begin{corollary}\label{cor:interpolation-degree}
The projective closure of $C_{n,d}$ in $\mathbb P^n$ has degree
\[
 \deg C_{n,d}=S(d+e-1,d).
\]
For general distinct $\alpha_1,\ldots,\alpha_{d+e-1}$ and general values
$\beta_1,\ldots,\beta_{d+e-1}$, there are exactly $S(d+e-1,d)$ such
polynomials satisfying $F(\alpha_i)=\beta_i$. Thus the incidence projection
$\pi$ of \cite[\S3]{vzGMat2024} has degree $\deg\pi=\deg C_{n,d}$.
\end{corollary}

\begin{proof}
Allow the coefficient of $t^{e-1}$ in $H$ and include the constant power
$H^0=1$. The same construction uses
\[
 \mathcal O\oplus\mathcal O(-1)\oplus\cdots\oplus\mathcal O(-d)
 \quad\hbox{over }\mathbb P^{e-1}.
\]
Its projectivization maps birationally onto the closure of $C_{n,d}$.
The Segre-class calculation gives $S(d+e-1,d)$.
Put $D=d+e-1$. The boundary of this closure in the hyperplane where the
leading coefficient vanishes has dimension at most $D-1$. A nonzero
polynomial has finitely many roots, so the incidence of a boundary point
with $D$ prescribed distinct zeros has dimension at most $D-1$.
Consequently, for general $(\alpha_i)$ the linear map
$[F]\mapsto[[t^n]F:F(\alpha_1):\cdots:F(\alpha_D)]$
has center disjoint from the closure. This linear projection to $\mathbb P^D$
is finite, pulls back $\mathcal O(1)$ to $\mathcal O(1)$, and has degree
$\deg C_{n,d}$. In characteristic zero its general fibers are reduced.
Uniqueness of the normalized decomposition identifies this degree with
that of the incidence projection.
\end{proof}

Thus the Stirling lower bound of \cite[Corollary~2.5 and Lemma~3.5]{vzGMat2024}
is the exact degree over $\C$. Centering removes one parameter, which
accounts for the different formula in \cref{th:scroll}.

\subsection{Composition chains}

For $\mathbf d=(d_1,\ldots,d_\nu)$, with $\nu\ge1$ and $d_i\ge2$, let
$\cC_{\mathbf d}$ consist of $F_1\circ\cdots\circ F_\nu$, where the factors
are monic, vanish at zero, and $F_\nu$ is centered. The standard uniqueness
and recovery results \cite{vzGMat2017,Ziegler} give
\[
 \cC_{\mathbf d}\simeq\A^{e_1},\qquad
 e_i=\sum_{j=i}^{\nu}(d_j-1)-1.
\]
Repeated use of the isolated-branch argument gives generic iterated
wreath-product monodromy. The invariant factors are indexed by the level
at which two leaves separate, so their weighted degrees are
\begin{equation}\label{eq:chain-type}
 \left(\left(\prod_{j=2}^{\nu}d_j\right)(d_1-1),\ldots,
 d_\nu(d_{\nu-1}-1),d_\nu-1\right).
\end{equation}

\begin{corollary}[Chain degree]\label{cor:chain-degree}
In the centered coefficient embedding,
\[
 \deg\overline{\cC}_{\mathbf d}
 =\prod_{i=1}^{\nu-1}S(d_i+e_{i+1},d_i).
\]
\end{corollary}

\begin{proof}
Build the tower from the innermost factor. Start with
$X_\nu=\mathbb P^{d_\nu-2}$ and $L_\nu=\mathcal O(-1)$. Recursively set
\[
 X_i=\mathbb P\left(\bigoplus_{j=1}^{d_i}L_{i+1}^{\otimes j}\right),
 \qquad L_i=\mathcal O_{X_i}(-1).
\]
The powers of each nonconstant polynomial are linearly independent, so
the construction gives a morphism from $X_1$ to the coefficient space.
Normalized uniqueness makes it birational onto
$\overline{\cC}_{\mathbf d}$. If $\xi_i=c_1(L_i^*)$, then
$\dim X_i=e_i$ and the projective-bundle formula gives
\[
 \int_{X_i}\xi_i^{e_i}
 =S(d_i+e_{i+1},d_i)\int_{X_{i+1}}\xi_{i+1}^{e_{i+1}}.
\]
Iteration proves the formula, including the empty product for $\nu=1$.
\end{proof}

\begin{proposition}\label{prop:dickson-degree}
For a prime $p\ge5$, the closure $\overline{\mathcal X}_p\subset\mathbb P^{p-2}$
of the Dickson curve is a rational normal curve of degree $(p-1)/2$ in its
linear span.
\end{proposition}

\begin{proof}
By the explicit formula for $D_p(t,a)$, the coefficient of $t^{p-2i}$ is a
nonzero multiple of $a^i$ for $0\le i\le(p-1)/2$, and all other coefficients
vanish. After homogenizing $a=a_1/a_0$, the parameterization
$a\mapsto D_p(t,a)$ is therefore given, in suitable coordinates, by the
monomials $a_0^{(p-1)/2-i}a_1^{i}$, $0\le i\le(p-1)/2$. Hence its image is a
rational normal curve of degree $(p-1)/2$ in a linear subspace of dimension
$(p-1)/2$.
\end{proof}

\subsection{Consequences and examples}

\begin{corollary}\label{cor:scroll-special}
We have $\deg\overline{\cD}_{r,2}=1$, $\deg\overline{\cD}_{2,s}=2^{s-1}-1$ and
$\deg\overline{\cD}_{r,3}=\binom{r+1}{2}$.
\end{corollary}

\begin{proof}
These are the Stirling numbers $\Stir{r}{r}$, $\Stir{s}{2}$ and
$\Stir{r+1}{r}$. The first equality can also be seen directly: for $s=2$ the
right composition factor is $t^2$, so $\cD_{r,2}$ is a linear subspace.
\end{proof}

\begin{example}\label{ex:scroll-two-three}
For $(r,s)=(2,3)$ one has
\[
 (t^3+ct)^2+b(t^3+ct)
 =t^6+2ct^4+bt^3+c^2t^2+bct.
\]
These polynomials form an affine plane embedded nonlinearly in
$\widetilde{\mathcal P}_6\simeq\A^4$; its projective closure is a cubic scroll surface in $\mathbb P^4$, in accordance with $\Stir{3}{2}=3$.
\end{example}

\subsection{The rational composition locus}\label{subsec:rational-scroll}

Let $k=rs$ with $r,s\ge2$, let $\cS$ be the tautological subbundle of rank two
on $\Gr(2,V_s)$, and put $\cW=\Sym^r\cS$, a bundle of rank $r+1$. Note that
$\dim V_s=s+1$, so $\Gr(2,V_s)\simeq\Gr(2,s+1)$ has dimension $2(s-1)$. By
the proof of \cref{prop:rational-composition}, the multiplication map
$\cW\to V_k\otimes\mathcal O$ is fiberwise injective, also at pencils with base
points. Hence it induces a morphism
\[
 \Psi\colon\ \Gr_2(\cW)\longrightarrow\Gr(2,V_k).
\]

\begin{theorem}[Rational composition degree]\label{th:rational-scroll}
Let $\mathcal C_s\subset\Gr(2,V_k)$ be the locus of base-point-free pencils of
degree $k$ that admit a right composition factor of degree $s$. The morphism
$\Psi$ is birational onto its image, and
\[
 \overline{\mathcal C_s}=\Psi(\Gr_2(\cW)),\qquad
 \dim\overline{\mathcal C_s}=2(r+s-2),\qquad
 \operatorname{codim}\overline{\mathcal C_s}=2(r-1)(s-1).
\]
Put $d=2(r+s-2)$, let $x,y$ be formal Chern roots of $\cS^*$, and define
\[
 w_i=(r-i)x+iy\quad(0\le i\le r),
\]
and
\[
 P_{r,s}(x,y)=\sum_{0\le i<j\le r}
 \frac{(w_i+w_j)^d}
 {\prod_{l\ne i,j}(w_i-w_l)(w_j-w_l)}.
\]
Then $P_{r,s}$ is a symmetric polynomial of degree $2s-2$, and the Pl\"ucker
degree of $\overline{\mathcal C_s}$ is
\begin{equation}\label{eq:rational-degree}
 \deg\overline{\mathcal C_s}
 =\int_{\Gr_2(\cW)}c_1(\det\cT^*)^d
 =[x^s y^{s-1}]\,\bigl((x-y)P_{r,s}(x,y)\bigr),
\end{equation}
where $\cT\subset\cW$ is the tautological subbundle of rank two on
$\Gr_2(\cW)$, that is, the pullback under $\Psi$ of the tautological subbundle
on $\Gr(2,V_k)$.
\end{theorem}

\begin{proof}
By \cref{prop:rational-composition}, the preimage under $\Psi$ of the
base-point-free locus of $\Gr(2,V_k)$ maps onto $\mathcal C_s$. This preimage
is a nonempty open subset of the irreducible Grassmann bundle $\Gr_2(\cW)$, and
$\Psi$ is proper, so the full image is $\overline{\mathcal C_s}$.

For birationality, let $g,h$ be general rational maps of degrees $r$ and $s$.
The noncollision argument used in the proofs of 
\cref{cor:rational-strata,th:generic-composition}, together with \cref{lem:wreath}, 
gives $\Mon(g\circ h)=S_s\wr S_r$. By \cref{lem:blocks},
such a composite $f$ has a unique right composition factor of degree $s$ up to
postcomposition, that is, the pencil $\cS_0=\langle A,B\rangle\in\Gr(2,V_s)$
with $h=A/B$ is uniquely determined by $f$. A point of $\Psi^{-1}(\Lambda_f)$,
where $\Lambda_f$ is the pencil of $f$, consists of a pencil $\cS_1$ such that
$\Lambda_f\subset\Sym^r\cS_1$ together with the two-plane $\Lambda_f$ itself,
viewed inside $\Sym^r\cS_1$; for base-point-free $\Lambda_f$, the pencil $\cS_1$
then defines a right composition factor of $f$, so $\cS_1=\cS_0$. Hence the
fiber of $\Psi$ over a general point of the image is a single point, and
$\Psi$ is birational onto its image, since we are in characteristic zero. The
dimension of $\Gr_2(\cW)$ is the sum of the base dimension $2(s-1)$ and the
fiber dimension $2(r-1)$; the codimension follows from
$\dim\Gr(2,V_k)=2(k-1)$.

The Pl\"ucker line bundle pulls back to $\det\cT^*$, so by birationality the
degree is $\int_{\Gr_2(\cW)}c_1(\det\cT^*)^d$. The Gysin formula for Grassmann
bundles \cite[Ch.~4]{FultonPragacz} expresses the pushforward of
$c_1(\det\cT^*)^d$ to $\Gr(2,V_s)$ as the sum $P_{r,s}$, where
$w_0,\dots,w_r$ are the Chern roots of $\cW^*=\Sym^r\cS^*$; the apparent
denominators cancel, because the pushforward of a polynomial in the Chern roots
of $\cT^*$ is a polynomial in the Chern classes of $\cW^*$. Finally, on
$\Gr(2,V_s)\simeq\Gr(2,s+1)$, the integral of a symmetric polynomial $F(x,y)$
of degree $2s-2$ in the Chern roots of $\cS^*$ equals
$[x^sy^{s-1}]\bigl((x-y)F(x,y)\bigr)$, as one checks on the basis of Schur
polynomials. This proves \eqref{eq:rational-degree}.
\end{proof}

\begin{remark}[Recovering the polynomial scroll]\label{rem:restriction}
Let $(t_0:t_1)$ be homogeneous coordinates on $\Pone$ with $t=t_1/t_0$, so that
the form $L=t_0^s\in V_s$ is the homogenization of the constant polynomial $1$,
and a polynomial map of degree $s$ corresponds to a pencil containing $L$. The
pencils containing $L$ form $\mathbb P(V_s/\langle L\rangle)\simeq\mathbb P^{s-1}$,
and the condition that the second generator be centered defines a hyperplane
$\mathbb P^{s-2}$ in it; a point $c$ of this hyperplane corresponds to the pencil spanned by $L$ and the homogenization of $\widehat H_c$ in $(t_0:t_1)$. On this subspace,
\[
 \cS|_{\mathbb P^{s-2}}=\mathcal O\oplus\mathcal O(-1),\qquad
 \Sym^r\cS|_{\mathbb P^{s-2}}=\mathcal O\oplus \cE,
\]
where the trivial summand is spanned by $L^r=t_0^{k}$. Requiring, in addition,
that the output two-plane contain $L^r$, that is, that the composite be a
polynomial, singles out the two-planes $\langle L^r,F\rangle$ with $F$ a line in
$\cE$; these form $\mathbb P(\cE)$, and we recover exactly \cref{th:scroll}.
Thus both the centering hyperplane and the constant line $L^r$ are needed to
pass from the rational Grassmann-bundle model to the polynomial one.
\end{remark}

The first values of \eqref{eq:rational-degree} are listed below. For
comparison, the centered polynomial model has degrees $\Stir{s}{2}=1,3,7,15$ for
$r=2$ and $s=2,3,4,5$.
\medskip

\begin{center}
\small
\renewcommand{\arraystretch}{1.15}
\begin{tabular}{@{}c|cccc@{}}
\toprule
$\deg\overline{\mathcal C_s}$&$s=2$&$s=3$&$s=4$&$s=5$\\
\midrule
$r=2$&$21$&$370$&$6\,125$&$98\,406$\\
$r=3$&$240$&$18\,060$&$1\,099\,700$&$59\,295\,852$\\
$r=4$&$2\,030$&$444\,150$&$71\,640\,800$&$9\,604\,985\,390$\\
$r=5$&$14\,490$&$7\,364\,280$&$2\,590\,487\,900$&$723\,573\,342\,492$\\
\bottomrule
\end{tabular}
\end{center}
\medskip

\begin{corollary}[Prescribed critical points]\label{cor:critical-points}
Let $D=2(r+s-2)$. The degree in \eqref{eq:rational-degree} counts degree-$rs$
rational maps with a right composition factor of degree $s$ and with $D$
prescribed general distinct critical points, modulo changes of the target
coordinate and counted with intersection multiplicity.
\end{corollary}

\begin{proof}
The Wronski map sends a pencil $\langle A,B\rangle$ to the
homogenization of $A'B-AB'$ in degree $2rs-2$. Its coordinates are linear in
the Pl\"ucker coordinates, and the form is nonzero for every pencil in
characteristic zero, including at pencils with base points. Thus prescribing
one critical point is a Pl\"ucker hyperplane condition; compare \cite{ErGab}.
Consider the incidence of a pencil in $\overline{\mathcal C_s}$ with an
ordered $D$-tuple of distinct zeros of its Wronski form. The projection to
the pencil has finite fibers, so this incidence has dimension at most $D$.
Hence a general $D$-tuple cuts out at most finitely many pencils.
The boundary has dimension at most $D-1$, and the same argument shows that
a general tuple avoids it. The resulting proper intersection of $D$
Pl\"ucker hyperplanes has intersection number
$\deg\overline{\mathcal C_s}$, as required.
\end{proof}

\begin{corollary}[Quadratic right composition factors]\label{cor:rational-quadratic}
For $s=2$,
\[
 \deg\overline{\mathcal C_2}
 =\frac{r(2r^2-r+1)}4\binom{2r}{r}.
\]
\end{corollary}

\begin{proof}
The base is $\Gr(2,3)=\mathbb P^2$. Write $h$ for its hyperplane class,
$a=r(r+1)/2$, and $\mathcal W=\Sym^r\cS^*$. Then
\[
 c_1(\mathcal W)=ah,\qquad c_2(\mathcal W)=\tfrac12a(a+1)h^2.
\]
In the Schur expansion of $(t_1+t_2)^{2r}$ on the Grassmann bundle, only
the partitions $(r+1,r-1)$ and $(r,r)$ survive pushforward. Their
coefficients are respectively
$\frac3{r+2}\binom{2r}{r-1}$ and
$\frac1{r+1}\binom{2r}{r}$. The pushforwards are
$c_1(\mathcal W)^2-c_2(\mathcal W)$ and $c_2(\mathcal W)$. Substitution and simplification give the
displayed formula. This also independently checks the column $s=2$.
\end{proof}

\section{Intersections of polynomial composition loci}\label{sec:ritt}

We use Ritt's theory here as background, rather than reprove the collision
classification. For a degree-$k$ polynomial, the branch cycle at infinity
is a $k$-cycle. Every monodromy block system is therefore a coset partition
of $\Z/k\Z$. There is at most one such system of each block size, and
the available sizes are closed under greatest common divisors and least
common multiples; see \cite[Corollary~2.6]{ZieveMueller}.

\begin{proposition}[Classical intersection reduction]\label{th:ritt-intersections}
Suppose $rs=r's'=k$ and $(r,s)\ne(r',s')$. Set $a=\gcd(r,r')$ and
$b=\gcd(s,s')$. There are coprime $m,n$ with
$r=am$, $r'=an$, $s=bn$, and $s'=bm$.
A polynomial lies in $\cD_{r,s}\cap\cD_{r',s'}$ precisely when, after
normalization, it has the form $A\circ Q\circ B$, with
$\deg A=a$, $\deg B=b$, and $Q$ admitting decompositions of patterns
$(m,n)$ and $(n,m)$.
\end{proposition}

This is the Ritt--Engstr\"om reduction; see
\cite{Engstrom,Tortrat} and \cite[Lemma~2.8 and Corollary~2.18]{ZieveMueller}.
Call the middle polynomial $Q$ in $A\circ Q\circ B$ the \emph{core}.
When $m,n>1$, Ritt's second theorem gives monic, centered representatives
$Q_0$ with $Q_0(0)=0$ in the affine equivalence class of $Q$. In fact, if
$Q=L_1\circ Q_0\circ L_2$ for affine automorphisms $L_1,L_2$, then
$A\circ Q\circ B=(A\circ L_1)\circ Q_0\circ(L_2\circ B)$, and the
degrees of $A$ and $B$ are preserved. With $\ell=\min(m,n)$ and $M=\max(m,n)=d\ell+e$,
the power and Dickson representatives $Q_0$ are
\[
 t^{e\ell}W(t^\ell)^\ell,\qquad
 D_{mn}(t,c)-D_{mn}(0,c),
\]
where $W$ is monic of degree $d$, $1\le e<\ell$, and $\gcd(e,\ell)=1$.
Precise normal forms and their overlaps are given in
\cite[Theorem~3.9]{vonzurGathen}.

For comparable right degrees, say $s'=qs$, the reduction gives
\begin{equation}\label{eq:nested}
 \cD_{qr',s}\cap\cD_{r',qs}=\cC_{(r',q,s)}.
\end{equation}
This is a smooth affine space of dimension $r'+q+s-4$. Its generic
invariant type is $(sq(r'-1),s(q-1),s-1)$, and \cref{cor:chain-degree} gives
the additional enumerative consequence
\[
 \deg\overline{\cC}_{(r',q,s)}
 =S(q+s-2,q)S(r'+q+s-3,r').
\]

\begin{corollary}[Singular locus]\label{cor:singular}
For composite $k$, with the reduced structure,
\[
 \Sing(\widetilde{\mathcal R}_k)
 =\bigcup_{\substack{rs=r's'=k\\(r,s)\ne(r',s')}}
   \cD_{r,s}\cap\cD_{r',s'}.
\]
For $k\ge3$, the locus is smooth precisely when $k$ is prime or the square
of a prime; for $k=3$ it is empty.
\end{corollary}

\begin{proof}
Each component is smooth. At a point belonging to two distinct components
the local ring has at least two minimal primes and is not regular.
All composition components contain $t^k$, and there is only one of them
precisely when $k$ is a prime square. In prime degree the Dickson component
is an affine line.
\end{proof}

We retain one explicit collision because it displays the invariant
factorization on an intersection. Its parameterization is a special case
of the classical power alternative.

\begin{proposition}\label{prop:swapped}
For $m=2d+1\ge3$,
\[
 \cD_{m,2}\cap\cD_{2,m}
 =\{t^2W(t^2)^2:W\text{ monic of degree }d\}\simeq\A^d.
\]
Writing $H=tW(t^2)$, the difference quotient factors as
\[
 B_{H^2}(z,w)=(z+w)
 \frac{H(z)-H(w)}{z-w}\frac{H(z)+H(w)}{z+w}.
\]
\end{proposition}

\begin{proof}
The first locus consists of even polynomials. A polynomial in both loci
can be written as $P=H^2+bH$, with $H$ the normalized right composition
factor of degree $m$. Then $K=H+b/2$ satisfies $K(-t)^2=K(t)^2$.
Its odd degree forces $K$ to be odd, whence $b=0$ and $H=tW(t^2)$.
The converse and the displayed factorization are immediate. Polynomial
recovery of $W$ gives the affine-space assertion.
\end{proof}

For $k=10$, put $P_{a,b}=t^2(t^4+at^2+b)^2$. Then
\begin{equation}\label{eq:k10}
\begin{split}
 \overline B_{P_{a,b}}(u,v)=u\,&
 (u^4-3u^2v+v^2+a(u^2-v)+b)\\
 &\cdot(u^4-5u^2v+5v^2+a(u^2-3v)+b).
\end{split}
\end{equation}
The generic type is $(4,4,1)$. The discriminants of the two quadratics in
$v$ are
$5u^4+2au^2+a^2-4b$ and $5u^4+10au^2+9a^2-20b$.
Each is a square in $\C[u]$ exactly when $b=a^2/5$. Thus both quartics
split precisely on this curve, where the type becomes $(2,2,2,2,1)$.
The substitution $a=-5c$, $b=5c^2$ identifies it with the Dickson family.

\section{Absolute and invariant irreducibility}\label{sec:prime}

We first recall the classical prime-degree classification, retaining only
the information needed for the factorization question. Burnside's theorem
reduces the non-$2$-transitive case to affine groups, and Ritt's genus-zero
analysis gives the following list; see \cite{Ritt3},
\cite[\S7]{Burda}, and \cite[\S6]{GuMuSa}.

\begin{proposition}[Classical prime-degree list]\label{th:prime-classification}
Let $f$ have prime degree $p\ge3$, with monodromy not $2$-transitive.
Then $G=\mathbb F_p\rtimes M$, where $M\le\mathbb F_p^*$ has order
$m\in\{1,2,3,4,6\}$, $m\mid p-1$, and $m<p-1$. The possibilities are
\[
\begin{array}{c|c|l}
 m&\text{signature of the Galois closure}&\text{map}\\\hline
 1&(p,p)&\text{power map}\\
 2&(p,2,2)&\text{Dickson map}\\
 2&(2,2,2,2)&\text{isogeny modulo }\mu_2\\
 3&(3,3,3)&\text{isogeny modulo }\mu_3\ (j=0)\\
 4&(2,4,4)&\text{isogeny modulo }\mu_4\ (j=1728)\\
 6&(2,3,6)&\text{isogeny modulo }\mu_6\ (j=0).
\end{array}
\]
The isogenies have degree $p$ and kernels invariant under $\mu_m$.
The residual curve has $(p-1)/m$ irreducible components, each of projection
degrees $(m,m)$; they are symmetric exactly when $m$ is even.
\end{proposition}

Identify a generic fiber with $\mathbb F_p$, on which $G$ acts by affine
maps. For an ordered pair of distinct points, let $x,y\in\mathbb F_p$ be
their labels. Translations preserve $y-x$, and multipliers in $M$ act on
it by multiplication, so the nontrivial orbitals are indexed by the cosets
$(y-x)M$ in $\mathbb F_p^*/M$. Transposition sends $(y-x)M$ to
$-(y-x)M$, so an orbital is self-paired precisely when $-1\in M$,
equivalently when $m$ is even.
The signatures follow from Riemann--Hurwitz: a translation has index
$p-1$, and a multiplier of order $d>1$ has index
$(p-1)(1-1/d)$. This is the calculation underlying the cited classical list.
In prime degree every nontrivial orbital is connected. Thus the list
also gives the irreducible finite-orbit correspondences with prime generic
orbit size and non-$2$-transitive quotient monodromy. The $2$-transitive
case gives the whole residual curve by \cref{th:dictionary}.
The only nonsymmetric cases are $m=1$ and $m=3$.

\begin{theorem}[The discrepancy in all degrees]\label{th:all-discrepancy}\label{th:kantor-sharp}
For a rational map $f\colon\Pone\to\Pone$ of degree at least two,
$B_f$ is reducible while $\overline B_f$ is irreducible if and only if,
up to M\"obius changes of source and target, one of the following holds:
\begin{enumerate}[{\rm (a)}]
\item $\deg f=3$ and $f(z)=z^3$;
\item $\deg f=7$ and $f$ is induced by a degree-seven isogeny of elliptic
curves with $j=0$, after quotienting by $\mu_3$; its kernel is
$\mu_3$-invariant, and $\Mon(f)=C_7\rtimes C_3$.
\end{enumerate}
In either case the absolute factors are the two transposed factors in
\eqref{eq:paired-factors}.
\end{theorem}

The maps in this statement are classical. In particular, the degree-seven
case belongs to Ritt's elliptic family; compare \cite[Theorem~1.4(b)]{GuMuSa}.
To pass from prime degree to arbitrary degree, we use the classical
classification of primitive affine genus-zero groups
\cite[Theorem~1.4]{Neubauer1992}; see also
\cite[Theorem~16 and Table~12]{MagaardShpectorovWang}.
This avoids an additional scalar-subgroup assumption in Kantor's theorem.

\begin{proof}
By \cref{prop:monodromy}, $G=\Mon(f)$ is $2$-homogeneous but not
$2$-transitive. It is primitive and has odd order: a nontrivial involution
would interchange a pair, making the single orbit on unordered pairs
a single orbit on ordered pairs. Degree two is impossible, and the degree-three case is immediate.
In degrees at least four, Kantor's theorem gives
$G\le\operatorname{A\Gamma L}_1(\mathbb F_q)$, where
$q=p^e\equiv3\pmod4$ \cite[Theorem~1]{Kantor}.

The group contains the full translation group
$V\simeq(\mathbb F_q,+)$. Indeed, $V\cap G$ is nontrivial since
$q\mid |G|$ but $q\nmid(q-1)e$; it is normal in the primitive group $G$,
so is transitive, and a transitive subgroup of the regular group $V$ is $V$.
Thus $G$ is a primitive affine genus-zero group with regular elementary
abelian socle of order $p^e$.

Neubauer's degree bounds \cite[Theorem~1.4]{Neubauer1992}, restated in
\cite[Theorem~16]{MagaardShpectorovWang}, give the alternatives
\[
 \begin{aligned}
 &1\le e\le2;\\
 &p=2,\ 2\le e\le8;\qquad p=3,\ 2\le e\le4;\\
 &p=5,\ 2\le e\le3;\qquad p\in\{7,11\},\ e=2.
 \end{aligned}
\]
Since $p^e\equiv3\pmod4$, $p$ and $e$ are odd. The only nonprime
possibility in this list is $q=27$. The six degree-$27$ groups in
\cite[Table~12]{MagaardShpectorovWang} have point stabilizers
$A_4$, $A_4\times C_2$, $S_4$, $S_4\times C_2$,
$\mathrm{SL}_3(3)$, or $\mathrm{GL}_3(3)$.
Each has even order, contradicting the odd order of $G$.

It follows that $q=p$ is prime. Then $G=\mathbb F_p\rtimes M$, and
$2$-homogeneity without $2$-transitivity means
$M\cup(-M)=\mathbb F_p^*$ with $-1\notin M$.
Thus $m=|M|=(p-1)/2$ is odd. In \cref{th:prime-classification}, the only
odd possibilities are $m=1,3$, giving exactly $p=3,7$ and the stated maps.
Conversely, these two groups have a single orbit on unordered pairs and
two transposed nontrivial orbitals. The degree-seven isogenies exist
because the order-three automorphism has eigenlines on the $7$-torsion.
\end{proof}

For completeness, the elliptic row $m=2$ has a useful explicit equation.
On $E:y^2=x^3+ax+b$, with $4a^3+27b^2\ne0$, put
\begin{equation}\label{eq:elliptic-factor}
 \Phi_{x_0}(z,w)=(zw+x_0(z+w)-a)^2
                 -4(z+w+x_0)(x_0zw+b).
\end{equation}
This is the classical $x$-coordinate addition relation: if $x_0=x(T)$,
its roots in $w$ are $x(R\pm T)$ when $z=x(R)$.
One obtains it by eliminating the slope and intercept of the line
through three points on $E$. See \cite[Theorem~16.1]{BuchsVesGa} for the
two-valued addition law and its degenerations, and
\cite[\S4]{MirmanShukla} for the Poncelet counterpart.
For a subgroup $C=\{O,\pm T_1,\ldots,\pm T_{(p-1)/2}\}$ of odd prime
order, the isogeny quotient satisfies
\[
 B_f=c\prod_{j=1}^{(p-1)/2}\Phi_{x(T_j)},\qquad c\ne0.
\]
Each factor is the irreducible symmetric orbital of differences
$\{\pm T_j\}$ and has projection degrees $(2,2)$.
Finiteness of these translation orbits is the classical torsion
condition; compare \cite[Theorem~2, Part~B]{FriedPoncelet}.

\section{Small degrees}\label{sec:small}

We now specialize the general results to small degrees; no degree-by-degree
elimination is needed. In the following table, dimension and codimension are
measured in $\widetilde{\mathcal P}_k\simeq\A^{k-2}$, ``generic type'' is the
generic invariant factorization type on the given component, and the last
column gives the degree of its projective closure, computed by
\cref{th:scroll} and \cref{prop:dickson-degree}.
\medskip

\begin{center}
\small
\renewcommand{\arraystretch}{1.18}
\begin{tabular}{@{}ccccc@{}}
\toprule
$k$&component of $\widetilde{\mathcal R}_k$&dim/codim&generic type&degree\\
\midrule
3&empty&--&--&--\\
4&$\cD_{2,2}$&$1/1$&$(2,1)$&$1$\\
5&$\mathcal X_5$&$1/2$&$(2,2)$&$2$\\
6&$\cD_{3,2}$&$2/2$&$(4,1)$&$1$\\
 &$\cD_{2,3}$&$2/2$&$(3,2)$&$3$\\
7&$\mathcal X_7$&$1/4$&$(2,2,2)$&$3$\\
8&$\cD_{4,2}$&$3/3$&$(6,1)$&$1$\\
 &$\cD_{2,4}$&$3/3$&$(4,3)$&$7$\\
9&$\cD_{3,3}$&$3/4$&$(6,2)$&$6$\\
10&$\cD_{5,2}$&$4/4$&$(8,1)$&$1$\\
  &$\cD_{2,5}$&$4/4$&$(5,4)$&$15$\\
11&$\mathcal X_{11}$&$1/8$&$(2,2,2,2,2)$&$5$\\
12&$\cD_{6,2}$&$5/5$&$(10,1)$&$1$\\
  &$\cD_{4,3}$&$4/6$&$(9,2)$&$10$\\
  &$\cD_{3,4}$&$4/6$&$(8,3)$&$25$\\
  &$\cD_{2,6}$&$5/5$&$(6,5)$&$31$\\
\bottomrule
\end{tabular}
\end{center}
\medskip

All rows follow from \cref{cor:complete-locus,th:generic-composition}, together with 
\cref{th:scroll,prop:dickson-degree} for the degrees: in
prime degree the only component is the Dickson curve, and in composite degree
the components are the composition loci, with generic type $(rs-s,s-1)$ and
degree $\Stir{r+s-2}{r}$.

The intersections are obtained from \cref{th:ritt-intersections}.
In degrees $6$ and $10$, \cref{prop:swapped} gives the affine families
of dimensions one and two. In degree $8$, the intersection is
$\cC_{(2,2,2)}$, of dimension two and generic type $(4,2,1)$.
The comparable-right-degree intersections in degree $12$ are
$\cC_{(3,2,2)}$, $\cC_{(2,3,2)}$, and $\cC_{(2,2,3)}$;
each has dimension three, with type given by \eqref{eq:chain-type}.
The noncomparable intersection $\cD_{4,3}\cap\cD_{3,4}$ is the union
of the power curve $\{t^3(t^3+c)^3\}$ and the Dickson curve, meeting
only at $t^{12}$. In every composite degree the Dickson curve lies
in all composition components.

\section{Further questions}

The main open factorization problem is to describe the reducibility locus for
general rational maps, in particular on the strata where the marked pole fiber
contains several ramified points. The Grassmann-bundle model of
\cref{subsec:rational-scroll} is a natural starting point for computing the
classes of the rational composition loci in the Chow ring of $\Gr(2,V_k)$. For
the polynomial composition loci the degree is determined by \cref{th:scroll};
their singularities at infinity and the intersection multiplicities of
different components are separate questions.

The appendix gives translation models with polynomial orbit growth. A
classification of nonsymmetric correspondences of projection degrees $(2,2)$
with polynomial orbit growth lies beyond the results proved here.

\appendix
\section{Translation examples}\label{sec:growth}

This appendix records the connection with multivalued groups without
attempting a classification of correspondences with slow orbit growth.
The quotient construction is classical \cite{Buchstaber,BuchsVes}; growth
of multivalued dynamics was studied by Veselov
\cite{Veselov1991,Veselov1992}, and comparison with growth in the covering
group is developed in \cite[Theorem~4.3]{BardakKozZon}.
Write $B_j(x)$ for the points reachable in at most $j$ forward or backward
steps. We count distinct endpoints, not paths with multiplicity.

Let $G$ be $\Ga$, $\Gm$, or an elliptic curve, written additively; let
$\Sigma\le\Aut(G,0)$ be finite with rational quotient, and let
$q:G\to G/\Sigma\subseteq\Pone$ be the quotient map.
For a nonempty finite $\Sigma$-stable set $S\subset G$, set
\begin{equation}\label{eq:translation-corr}
 \Gamma=\overline{\{(q(g),q(g+s)):g\in G,\ s\in S\}}.
\end{equation}

\begin{proposition}\label{prop:growth-rank}
The components of $\Gamma$ correspond to the $\Sigma$-orbits in $S$.
Its projection degrees are $(|S|,|S|)$, and it is symmetric exactly when
$S=-S$. For every $x\in q(G)$,
\[
 |B_j(x)|\asymp j^{\rank\Lambda},\qquad \Lambda=\langle S\rangle,
\]
where rank means the rank of the finitely generated abelian group
$\Lambda$ modulo its torsion subgroup.
\end{proposition}

\begin{proof}
Two translation graphs have the same image under $q\times q$ exactly
when their steps are in the same $\Sigma$-orbit. For general $g$, the
values $q(g+s)$ are distinct: a collision with a nonidentity
$\sigma\in\Sigma$ requires $(1-\sigma)g$ to take a prescribed value,
and $1-\sigma$ has finite kernel. This proves the component and degree
assertions. Transposition replaces $S$ by $-S$.
Every path lifts to successive additions of elements of $S\cup(-S)$,
because $S$ is $\Sigma$-stable. Conversely, every such sequence projects
to a path. The map $q$ identifies at most $|\Sigma|$ points, so the
endpoint count lies between $1/|\Sigma|$ times the size of the word ball
in $\Lambda$ and that size itself. Such balls grow as
$j^{\rank\Lambda}$. For nonproper $G$, translation preserves the
boundary of its smooth completion; no added edge joins $q(G)$ to the
boundary.
\end{proof}

The familiar symmetric biquadratics are the elliptic relation
$x(R)\leftrightarrow x(R\pm T)$, with $T\notin E[2]$, and its
multiplicative and additive degenerations
\[
 z^2-(\lambda+\lambda^{-1})zw+w^2
       +a(\lambda-\lambda^{-1})^2=0,
 \qquad (w-z-1)^2=4z,
\]
where $a\ne0$ and $\lambda\ne\pm1$ in the first equation.
Up to a simultaneous M\"obius change, these exhaust the irreducible
symmetric curves of projection degrees $(2,2)$.
These are classical normal forms; compare
\cite[Theorem~2, equations~(19)--(21)]{ABS} and
\cite[Theorem~16.1]{BuchsVesGa}. One way to see the trichotomy is to
normalize the curve, whose genus is at most one. The product of the
first projection's deck involution and the coordinate-flip involution
is a translation in genus one and a semisimple or parabolic M\"obius
map in genus zero. The elliptic and multiplicative orbits are finite
exactly for torsion steps, as in \cite[Theorem~2, Part~B]{FriedPoncelet};
otherwise the orbits in the affine quotient grow linearly. The additive
case has linear growth, and the rational models have a singleton orbit
at infinity.

\begin{corollary}[Cyclotomic steps]\label{th:phi-growth}
For $d\ge1$ and $c\ne0$, the curve
$\Gamma_d=\overline{\{(g^d,(g+c)^d):g\in\C\}}$ is irreducible of
projection degrees $(d,d)$ and is symmetric exactly when $d$ is even.
The full orbit of $g^d\in\C$ is
\[
 \{(g+\lambda)^d:\lambda\in c\Z[\zeta_d]\},
\]
and its growth is $\asymp j^{\varphi(d)}$.
There is a nonconstant meromorphic first integral on $\C$ precisely
when $d\in\{1,2,3,4,6\}$. The point $\infty$ is a singleton orbit.
\end{corollary}

\begin{proof}
Use \cref{prop:growth-rank} with $G=\Ga$, $\Sigma=\mu_d$, and
$S=\mu_dc$. The parameterization extends to a morphism
$\Pone\to\Pone\times\Pone$ whose only point over the boundary is
$\infty\mapsto(\infty,\infty)$. Its image is irreducible and closed,
which also proves the assertion about $\infty$.
The step group $\Lambda=c\Z[\zeta_d]$ has rank $\varphi(d)$.
A discrete subgroup of $\C$ has rank at most two, and $\varphi(d)\le2$
holds exactly for $d\in\{1,2,3,4,6\}$. For these values $\Lambda$ is
cyclic or a lattice, so it is indeed discrete.
In the lattice cases, pull back a nonconstant meromorphic function on
the compact Riemann surface $(\C/\Lambda)/\mu_d$; in the cyclic cases
$d=1,2$, use $\exp(2\pi i g/c)$ and $\cos(2\pi g/c)$, respectively.
These functions are $\Lambda$-periodic and $\mu_d$-invariant, and hence
descend to meromorphic functions of $g^d$
\cite[Theorem~2.1 and Corollary~14.13]{Forster}.
Conversely, a meromorphic first integral lifts to a function on $\C$
with every element of $c\Z[\zeta_d]$ as a period. If the group is not
discrete, periods tend to zero, forcing the derivative to vanish away
from the poles. The function is then constant.
\end{proof}

The dynamics $\Gamma_d$ belongs to the Buchstaber--Novikov family studied
by Chirkov \cite{Chirkov}. The growth exponent above follows from the
rank of its step group; the general theory of finite sumsets gives
stronger eventual-polynomial statements for the covering group
\cite{Khov}. No general classification of nonsymmetric biquadratics
or of correspondences with polynomial orbit growth is claimed here.

\medskip
\noindent\textbf{Acknowledgements.}
The second author thanks A.~Veselov, I.~Sechin, and N.~Hemmingsson for useful
discussions.

\medskip
\noindent\textbf{Declaration of generative AI and AI-assisted technologies
in the manuscript preparation process.}
Claude Sonnet (Anthropic) was used to improve the original manuscript on the
polynomial case and to generate results for the rational case. GPT-based
tools (ChatGPT and Codex, OpenAI) were used for literature searches, to
suggest and check additional mathematical arguments and formulas, to revise
the exposition, and to assist with partial Lean formalizations. The authors
checked the results and take responsibility for the final content of this
article.

\end{document}